\documentclass[10pt]{article}

\usepackage{amssymb}
\usepackage{amsmath}
\input{liemacs.sty} 

\begin{document} 

%%%%%%%%%%%%%%%%%%%%%%%%%%

%\title{Holomorphically induced representations of root graded 
%Banach--Lie groups} 
\title{Borel--Weil Theory for Root Graded Banach--Lie Groups} 

\author
{Christoph M\"uller\footnote{Technische Universit\"at Darmstadt,
Schlossgartenstrasse 7, D-64289 Darmstadt, Deutschland,
cmueller@mathematik.tu-darmstadt.de. The first author gratefully 
acknowledges support by the DFG-Project 
``Geometrische Darstellungstheorie wurzelgraduierter 
Lie-Gruppen'', Ne 413/5-2.},
Karl-Hermann
Neeb\footnote{Technische Universit\"at Darmstadt,
Schlossgartenstrasse 7, D-64289 Darmstadt, Deutschland,
neeb@mathematik.tu-darmstadt.de},
Henrik Sepp\"anen\footnote{Technische Universit\"at Darmstadt,
Schlossgartenstrasse 7, D-64289 Darmstadt, Deutschland,
seppaenen@mathematik.tu-darmstadt.de. The third author was supported
by a post-doctoral fellowship from the Swedish Research Council.}}

%\address{Fachbereich Mathematik, AG AGF
%Technische Universit\"at Darmstadt
%Schlossgartenstrasse 7
%64289 Darmstadt}
%\email{seppaenen@mathematik.tu-darmstadt.de}
%\keywords{}
%\subjclass{}

\newcommand{\Ind}{\mathop{{\rm Ind}}\nolimits}
\newcommand{\Coind}{\mathop{{\rm Coind}}\nolimits}
\newcommand{\bF}{{\mathbb F}}
\newcommand{\cR}{{\mathcal R}}
\newcommand{\E}{{\mathbb E}}

\maketitle

\begin{abstract} 
In this paper we introduce (weakly) root graded Banach--Lie 
algebras and corresponding Lie groups as natural generalizations 
of group like $\GL_n(A)$ for a Banach algebra $A$ or groups like 
$C(X,K)$ of continuous maps of a compact space $X$ into a complex 
semisimple Lie group~$K$. 
We study holomorphic induction from holomorphic Banach 
representations of so-called parabolic subgroups $P$ 
to representations of $G$ on holomorphic sections 
of homogeneous vector bundles over $G/P$. One of our main results 
is an algebraic characterization 
of the space of sections which is used to show that this space 
actually carries a natural Banach structure, a result generalizing 
the finite dimensionality of spaces of sections of holomorphic bundles 
over compact complex manifolds. 
We also give a geometric realization of any irreducible holomorphic 
representation of a (weakly) root graded Banach--Lie group $G$ 
and show that all holomorphic functions on the spaces $G/P$ are constant.
\end{abstract}

\section*{Introduction} 

Let $\g_\Delta$ be the 
finite-dimensional semisimple complex Lie algebra with 
root system $\Delta$ and fix a Cartan subalgebra $\h$ of $\g_\Delta$. 
A complex Banach--Lie algebra $\g$ is said to be 
weakly $\Delta$-graded if it contains $\g_\Delta$ and 
decomposes as a direct sum 
$$ \g = \g_0 \oplus \bigoplus_{\alpha \in \cR} \g_\alpha $$
of simultaneous $\ad \h$-eigenspaces.  
It is called {\it root graded} if, in addition, $\Delta = \cR$ and $\g$ is 
topologically generated by the root spaces $\g_\alpha$. 

The systematic study of root graded Lie algebras (with irreducible $\Delta$) 
was initiated by 
Berman and Moody in \cite{BM92}, where they studied Lie algebras 
graded by  root systems of the  types $A$, $D$ and $E$. 
Corresponding results for non-simply laced root systems have been obtained by 
Benkart and Zelmanov in \cite{BZ96}. 
The classification of root graded Lie algebras in the algebraic 
context was completed 
by Allison, Benkart and Gao in \cite{ABG00}. 
The main point of the classification is to associate to a 
$\Delta$-graded Lie algebra $\g$ a coordinate algebra $\cA$ of a 
certain type depending on $\Delta$, and then to show that, up to 
central extensions, $\g$ is determined by its coordinate algebra. 
For type $A_1$, the coordinate algebras turn out to be unital 
Jordan algebras, for $A_2$ unital alternative algebras, and for 
$A_n$, $n > 2$, they are unital associative algebras. 
For types $D$ and $E$ they are commutative, so that $\g$ is a 
central extension of $\cA \otimes \g_\Delta$. 
Apart from  simple complex Lie algebras, the most well-known class of 
root graded Lie algebras are affine Kac--Moody algebras (\cite[Ch.~6]{Ka90}). 
For the case $\cR = BC_r$ and $\Delta$ of type $B$, $C$ or $D$, 
we refer to the memoir \cite{ABG02}. 
The classification scheme for root graded Lie algebras has been extended 
to the topological context of locally convex Lie algebras 
in \cite{Ne03} to cover many classes 
arising mathematical physics, 
operator theory and geometry.

This paper is the first in a series dedicated to various 
aspects of holomorphic representations of weakly 
root graded Banach--Lie groups $G$, 
i.e., groups whose Lie algebra $\g$ is weakly root graded. 
Since there is a natural notion of parabolic subgroups $P$ of a 
root graded Lie group $G$, and the corresponding homogeneous spaces 
$G/P$ carry complex manifold structures, it is a natural 
problem to understand the representations of $G$ in spaces of 
holomorphic sections of homogeneous holomorphic vector bundles 
$\bE_\rho = G \times_\rho E$ over $G/P$, defined by a holomorphic representation 
$\rho \: P \to \GL(E)$, where $E$ is a Banach space and 
$\rho$ is a morphism of Banach--Lie groups. 

The classical context for these problems is the 
Borel--Weil Theorem, where $G$ is a complex reductive Lie group 
and $G/P$ is a generalized flag manifold, hence in particular 
compact. In this case the space of holomorphic sections is always 
finite dimensional if $E$ is so. The Borel--Weil Theorem 
identifies those holomorphic characters $\rho \: P \to \C^\times = \GL(\C)$ 
for which this space is non-zero, and the natural $G$ representation on this 
space, which is irreducible. 

One of the main results of this paper is that for each 
holomorphic Banach representation $(\rho,E)$ of a parabolic 
subgroup $P$, the space of holomorphic sections of $\bE_\rho$, 
which we identify with the space 
$$ \cO_\rho(G,E) 
= \{ f \in {\cal O}(G,E) \: (\forall g \in G)(\forall p \in P)\ 
f(gp) = \rho(p)^{-1}f(g)\} $$
carries a natural Banach space structure for which the 
representation of $G$, defined by $(\pi(g)f)(x) = f(g^{-1}x)$ 
defines a morphism $\pi \: G \to \GL(\cO_\rho(G,E))$ of 
Banach--Lie groups (Theorem~\ref{thm:fin-weight}). 
Here the remarkable point is that on the much larger space 
$\cO(G,E)$ of all $E$-valued holomorphic functions, 
there is no natural locally convex structure for which the 
$G$-action is continuous, resp., holomorphic. A natural topology 
on this space is the compact open topology, i.e., the topology of 
uniform convergence on compact subsets of~$G$. With respect to this 
topology, $\cO(G,E)$ is a complete locally convex space 
(\cite[Thm.~III.11(c)]{Ne01}), but for infinite dimensional 
groups the action of $G$ is not continuous. However, the finite dimensional 
complex 
subgroup $G_\Delta$ corresponding to the subalgebra $\g_\Delta$ of $\g$ 
acts holomorphically on this space. An instructive example 
illustrating the problem is the action of a Banach space 
$X$, considered as an abelian Banach--Lie group, on the space 
$A$ of affine holomorphic functions $X \to \C$ by 
$\pi(x)(f)(y) = f(y-x)$. This action is continuous with respect to the 
compact open topology if and only if the evaluation map 
$X \times X' \to \C, (x,f) \mapsto f(x)$ is continuous, but this is 
equivalent to $X$ being finite dimensional because otherwise 
the bipolars of compact subsets of $X$ are never $0$-neighborhoods 
(cf.\ \cite[p.2]{KM97}, \cite[2.2.10-13]{Mue06}). 
Our strategy to find a better 
topology on the space $\cO_\rho(G,E)$ is to realize it, 
basically by Taylor expansion in~$\1$, as a set of linear 
maps $\alpha \: U(\g) \to E$ 
on the enveloping algebra $U(\g)$ of the Lie algebra $\g$ of $G$. 
These linear maps are always {\it continuous} in the 
sense that all $n$-linear maps  
$\alpha_n \: \g^n \to E, (x_1,\ldots, x_n) \mapsto 
\alpha(x_1 \cdots x_n)$ are continuous, and the space 
$\Hom(U(\g),E)_c$ of all continuous linear maps carries a natural 
Fr\'echet space structure with respect 
to which the natural $\g$ action is a continuous bilinear map  
(Section~\ref{sec:2}). The key observation is that in the 
coinduced representation $\Hom_\fp(U(\g),E)_c$ corresponding to the 
$\fp$-module $E$, the space of $\g_\Delta$-finite vectors 
decomposes into finitely many $\fh$-weight spaces. Based on this 
observation, we proceed to show that it is closed in the 
aforementioned Fr\'echet space and that it even is a Banach space. 
This provides the infinitesimal picture. 
The bridge to the $G$-action on $\cO_\rho(G,E)$ is developed 
in Section~\ref{sec:3}. It is based on an application of the 
general Peter--Weyl Theorem to the action of a compact 
real form of the semisimple complex Lie group $G_\Delta$ 
on the locally convex spaces 
$\cO_\rho(G,E)$, endowed with the compact open topology. 
If $G$ is $1$-connected and $P$ is connected, we show that the 
map 
$$ \Phi \: \cO_\rho(G,E) \to \Hom_\fp(U(\g),E)_c, \quad 
\Phi(f)(D) := (D_r f)(\1), $$
where $D_r$ is the right invariant differential operator on $G$ 
associated to $D \in U(\g)$, defines a bijection onto the Banach 
subspace of $\g_\Delta$-finite elements in the Fr\'echet space 
$\Hom_\fp(U(\g),E)_c$ (Theorem~\ref{thm:fin-weight}). 
This result has some immediate consequences that are derived in 
Section~\ref{sec:3}. One is that all holomorphic functions on 
$G/P$ are constant, although this manifold is far from being compact. 
Another one is that $\cO_\rho(G,E)$ is finite dimensional if 
$G$ and $E$ are. The proof is rather elementary and does not
rely on any deep theory of elliptic operators or sheaves.

In connected complex reductive Lie groups, the parabolic subgroups 
are always connected. This is no longer the case for root graded 
Lie groups, so that one has to understand the influence of the 
passage from a holomorphic representation $\rho \: P \to \GL(E)$ 
to the restriction $\rho_0 := \rho\res_{P_0}$ to its identity 
component. Clearly, $\cO_\rho(P,E) \subeq \cO_{\rho_0}(P,E)$, 
but it is also natural to start with a representation $\rho_0$ of $P_0$ 
and ask for extensions to $P$ for which $\cO_\rho(G,P)$ is non-zero. 
In Theorem~\ref{thm:extens} we give a complete answer to this question 
in the most important case, where $\rho_0$ is irreducible with 
$\End_{P_0}(E) = \C \1$. We show that if all representations 
$p.\rho_0$, defined by 
$(p.\rho_0)(x) := \rho_0(p^{-1}xp)$ are equivalent to $\rho_0$ 
(which is necessary for the existence of some extension $\rho$) 
and $\cO_{\rho_0}(G,E)$ is non-zero, then  there exists a unique 
extension $\rho$ to $P$ with $\cO_\rho(P,E) = \cO_{\rho_0}(G,P)$, 
and for all other extensions $\gamma$ of $\rho_0$, the space 
$\cO_\gamma(G,P)$ vanishes. We find this quite surprising. 
It reduces all questions on the representations 
$\cO_\rho(G,E)$ to the case where $P$ is connected and 
$G$ is simply connected, and then $\Phi$ maps it isomorphically 
to $\Hom_\fp(U(\g),E)_c^{[\g_\Delta]}$. 
In Section~\ref{sec:5} we show that for each connected parabolic 
subgroup $P$, each irreducible 
holomorphic representation of $G$ can be realized in 
some space $\cO_\rho(G,E)$ for an irreducible holomorphic representation 
$(\rho,E)$ of $P$. The difficult question that remains open at this point 
is a characterization of those irreducible holomorphic representations 
$(\rho,E)$ of $P$ for which $\cO_\rho(G,E)$ is non-zero, which implies 
the existence of a corresponding irreducible holomorphic $G$-representation 
sitting as a minimal non-zero submodule in $\cO_\rho(G,E)$. 

In an appendix we develop a quite general variant of  
Frobenius Reciprocity for representations of $G$ on 
locally convex spaces for which all orbit maps are holomorphic. 
This applies in particular to the representation on 
$\cO(G,E)$, endowed with the topology of pointwise convergence. 

We plan to continue this project in subsequent papers 
which address the special cases where 
$G$ is finite dimensional but not necessarily semisimple, such 
as $\GL_n(A)$ for a finite dimensional algebra $A$, and 
where $\dim E = 1$ and $G = \GL_n(A)$ for a Banach algebra~$A$, resp., 
$\g = A \otimes \g_\Delta$ for a commutative Banach algebra~$A$. 
The main result is a characterization of those homogeneous 
line bundles which admit global holomorphic sections. It 
turns out that these are generated by pullbacks of positive 
line bundles over compact flag manifolds $G/P$ with respect 
to natural mappings induced by characters of the algebra $A$.

Several results presented in the present paper grew out of 
predecessors from Ch.~M\"uller's thesis \cite{Mue06} which 
deals only with the scalar case $E =~\C$. 

\section{Root graded Banach--Lie groups}\label{sec:1}

In this section we introduce weakly root graded Lie algebras and their 
parabolic subalgebras. We also discuss parabolic subgroups of corresponding 
Banach--Lie groups and how they can be used to obtain triple coordinates 
of large open identity neighborhoods. 

\subsection*{Root graded Lie algebras and parabolic subalgebras} 

\begin{definition} \label{def:1.1} 
Let $\Delta$ be a finite reduced root system 
and $\g_\Delta$ be the corresponding 
finite dimensional complex semisimple Lie algebra. 
A complex Banach--Lie algebra $\g$ is said to be 
{\it weakly $\Delta$-graded} if the following conditions are satisfied:
\begin{description}
\item[\rm(R1)] $\g_\Delta$ is a Lie subalgebra of $\g$. 
\item[\rm(R2)] For some (and hence for each) Cartan subalgebra 
$\fh$ of $\g_\Delta$, $\g$ decomposes as a finite direct sum 
of $\fh$-eigenspaces 
$\g = \g_0 \oplus \bigoplus_{\alpha \in \cR} \g_\alpha$, where
$\g_\alpha = \{x \in \g \: (\forall h \in \h)\, [h,x] =
\alpha(h)x\}$. 
\end{description}

The subalgebra $\g_\Delta$ of $\g$ is called a {\it grading subalgebra}. 
We say that $\g$ is {\it (weakly) root graded} if $\g$ is (weakly) 
$\Delta$-graded for some $\Delta$. 
If, in addition, $\cR = \Delta$ and 
$\sum_{\alpha \in \Delta} 
[\g_\alpha, \g_{-\alpha}]$ is dense in $\g_0$, 
then $\g$ is said to be {\it $\Delta$-graded}. 

For $\alpha \in \Delta$, the 
unique element $\check \alpha \in 
[\g_{\Delta,\alpha}, \g_{\Delta,-\alpha}]$, 
satisfying $\alpha(\check \alpha) = 2$ is called the {\it 
coroot} corresponding to $\alpha$. From the representation theory of 
$\fsl_2(\C)$, it follows that $\cR(\check \alpha) \subeq \Z$ for 
each coroot, and it is well known that $\fh = \Spann \check\Delta$. 
\end{definition}

\begin{examples} (a) Let $\Delta$ be a reduced finite root system and 
$\g_\Delta$ be the corresponding semisimple complex Lie algebra. 
If $A$ is a commutative unital Banach algebra, then 
$\g := A \otimes \g_\Delta$ is a $\Delta$-graded 
Banach--Lie algebra with respect to the 
bracket 
$$ [a \otimes x, a' \otimes x'] := aa' \otimes [x,x']. $$
The embedding $\g_\Delta \into \g$ is given by $x \mapsto \1 \otimes x$. 

(b) If $A$ is a unital Banach algebra, then 
the $(n \times n)$-matrix algebra $M_n(A) \cong A \otimes M_n(\C)$ 
also is a Banach algebra. We write $\gl_n(A)$ for this algebra, 
endowed with the commutator bracket. 
Then $\gl_n(A)$ is a weakly $A_{n-1}$-graded Lie algebra 
with grading subalgebra $\g_\Delta = \1 \otimes \fsl_n(\C)$. 

(c) Let $\g = \g_1 \oplus \g_0 \oplus \g_{-1}$ be a 
$3$-graded Banach--Lie algebra for which there exist 
elements $e \in \g_1$ and $f \in \g_{-1}$ such that 
$h := [e,f]$ satisfies 
$$ \g_{\pm 1} = \{ x \in \g \: [h,x] = \pm 2 x\}. $$ 
Then $\g$ is weakly $A_1$-graded with grading subalgebra 
$\g_\Delta = \Spann \{ e,f,h\}$ (cf.\ \cite{BN04}). 

If, more generally, $\g = \sum_{j = -n}^n \g_j$ 
is $(2n+1)$-graded and there exist 
$e \in \g_2$ and $f \in \g_{-2}$ such that 
$h := [e,f]$ satisfies 
$$ \g_j = \{ x \in \g \: [h,x] = j x\}, $$ 
then $\g$ is weakly $A_1$-graded. 
%Note that the representation 
%theory of $\sL_2(\C)$ implies that the only possible eigenvalues 
%of $\ad \fh$ on $\g$ are integral. 

(d)  Assume that $\g$ is weakly $\Delta$-graded, and that $V$ is a 
$\g$-module which decomposes into a direct sum of finitely many 
weight spaces under $\h$.
Then the trivial abelian extension $V \rtimes \g$ of $\g$
by $V$ defined as 
$$[(v_1,x_1),(v_2,x_2)]:=(x_1.v_2+x_2.v_1,[x_1,x_2]), 
\quad v_1, v_2 \in V, x_1, x_2 \in \g$$
is also weakly $\Delta$-graded. 

\end{examples}

\begin{rem} \label{rem:centext} (Central extensions) 
(a) If $\g$ is (weakly) $\Delta$-graded, then 
its center $\z(\g)$ is contained in $\g_0$ and intersects 
$\g_\Delta$ trivially, so that the Banach--Lie algebra $\g/\z(\g) 
\cong \ad \g$ is also (weakly) $\Delta$-graded. 

(b) If $q \: \hat\g \to \g$ is a central extension of the 
weakly $\Delta$-graded Banach--Lie algebra $\g$ with kernel $\fz$, 
then the adjoint 
action of $\g$ on itself lifts to a natural action on $\hat\g$ for which 
$q$ is equivariant. Since the central extension 
$q^{-1}(\g_\Delta)$ of $\g_\Delta$ splits, we may consider 
$\g_\Delta$, and therefore also $\fh$, as a subalgebra of $\hat\g$. 
Then the subspaces $q^{-1}(\g_\alpha)$, $\alpha \in \cR$, are 
$\fh$-invariant and contain $\fz \subeq \hat\g_0$. 
If $x \in \fh$ satisfies $\alpha(x) = 1$, then 
$(\ad x- \1)(q^{-1}(\g_\alpha)) \subeq \ker q = \fz$ implies that 
$\ad x (\ad x - \1)$ vanishes on $q^{-1}(\g_\alpha)$, hence that 
$\ad x$ is diagonalizable on this space with eigenvalues $0$ and $1$. 
Clearly $\z \subeq \ker \ad x$ and 
$[x, q^{-1}(\g_\alpha)]$ maps surjectively onto 
$\g_\alpha$, which leads to the direct decomposition 
$$ q^{-1}(\g_\alpha) = \z \oplus [x,q^{-1}(\g_\alpha)]. $$
Since $\fh$ is abelian, $[x,q^{-1}(\g_\alpha)]$ is $\fh$-invariant, 
and we derive that 
$$ \hat\g_\alpha = \{ y \in \hat\g \: (\forall h \in \fh)\, 
[h,y] = \alpha(h)y \} = [x,q^{-1}(\g_\alpha)]. $$

This proves that $\ad \fh$ is diagonalizable on 
the subspace $\fz + \sum_{\alpha \in \cR} \hat\g_\alpha \subeq \hat\g$. 
If we assume, in addition, that $\hat\g_0 := q^{-1}(\g_0)$ is centralized 
by $\fh$ (which is automatic if $\hat\g$ is generated by  the subspaces 
$\hat\g_\alpha$), then $\hat\g$ is also weakly $\Delta$-graded. 
\end{rem}

For more details on other classes of examples and a discussion of 
topological issues related to root graded Lie algebras, we refer 
to \cite{Ne03}. 

\begin{defn} \label{def:parasys} A subset $\Sigma \subeq \cR$ is called a 
{\it parabolic systems} if 
there exists an element 
$x \in \fh$ with 
\begin{equation}
  \label{eq:pardef}
\Sigma = \{ \alpha \in \cR \: \alpha(x) \geq 0\}. 
\end{equation}
For a parabolic system $\Sigma$, we put 
$$ \Sigma ^+ := \Sigma \setminus -\Sigma, \quad 
\Sigma^0 := \Sigma \cap  -\Sigma \quad \mbox{ and } \quad 
\Sigma^- = \cR \setminus \Sigma. $$
For each parabolic system $\Sigma$, 
$\fp_\Sigma := \g_0 + \sum_{\alpha \in \Sigma} \g_\alpha$
is a subalgebra of $\g$, called 
the corresponding {\it parabolic subalgebra}. 

A parabolic system $\Sigma$ is called a {\it positive system} if 
$\Sigma^0 = \eset$. Since $\cR$ is real-valued on 
$\fh_\R := \Spann_\R \check\Delta$, the set of all elements 
$y \in \fh_\R$ with $\cR(y) \subeq \R^\times$ is dense, and 
choosing $y$ sufficiently close to the element $x$ from above, it follows that 
each parabolic system $\Sigma$ contains 
a positive systems $\cR^+$. 
%
%If $\cR^+ \subeq \Sigma$ is a positive system and 
%$\Pi \subeq \Delta^+ := \Delta \cap \cR^+$ is the set of simple roots, then 
%$\Pi$ is a basis of $\h^*$ and there exists a subset 
%$\Pi_\Sigma \subeq \Pi$ with 
%$$ \Sigma \cap \Delta = \Delta^+ \cup \Spann_\N(-\Pi_\Sigma). $$
\end{defn}

\begin{remark} \label{rem:para-pos} 
Since $\cR$ is finite, there exists for each parabolic system 
$\Sigma$ an element $x_\Sigma \in \fh$ with 
$\Sigma^+(x_\Sigma) \geq 1$ and 
$\Sigma^-(x_\Sigma) \leq -1$. 
\end{remark}

\subsection*{Root graded Lie groups and parabolic subgroups} 

\begin{definition} A Banach--Lie group $G$ is said to be 
{\it (weakly) root graded} if its Lie algebra $\L(G)$ is (weakly) 
root graded. 
\end{definition} 

\begin{rem} (a) Not every Banach--Lie algebra $\g$ is integrable 
in the sense that there exists a Banach--Lie group $G$ with 
Lie algebra $\g$. Since $\g = C^1(\bS^1,\fsl_2(\C))$ 
is $A_1$-graded with grading algebra $\fsl_2(\C)$, 
and this property is inherited by the central extension 
defined by the $\fsl_2(\C)$-invariant cocycle 
$$ \omega(\xi,\eta) = \int_0^{2\pi} \tr(\xi(t) \eta'(t))\, dt $$ 
(cf.\ Remark~\ref{rem:centext}), 
the discussion in \cite{EK64} implies the existence of non-integrable 
root graded Banach--Lie algebras. 

(b) If $A$ is a commutative unital Banach algebra, then 
any inclusion $\g_\Delta 
\subeq \gl_n(\C)$ yields an inclusion 
$\g \into A \otimes M_n(\C) \cong M_n(A)$ of Banach--Lie algebras, 
which implies the existence of a Lie group $G$ with Lie algebra $\g$ 
(\cite{Ms62}, \cite[Thm.~IV.4.9]{Ne06}). 

(c) For similar reasons, for each Banach--Lie algebra $\g$, the 
quotient $\g/\z(\g)$ is integrable because it embeds continuously 
into the integrable Lie algebra $\der(\g) = \L(\Aut(\g))$ 
(\cite{EK64}, \cite{Ne06}). 
\end{rem}

In the following, $G$ is a weakly $\Delta$-graded Banach--Lie group. 
If $\fp$ is a parabolic subalgebra, we write 
$\Sigma$ for the corresponding parabolic system with 
$\fp= \g_0 + \sum_{\alpha \in \Sigma} \g_\alpha.$
Then we have a semidirect decomposition 
$$ \fp = \fu \rtimes \fl, \quad \mbox{ where } \quad 
\fl= \g_0 + \sum_{\alpha \in \Sigma^0} \g_\alpha \quad \mbox{ and } \quad 
\fu= \sum_{\alpha \in \Sigma^+} \g_\alpha. $$
That $\fl$ is a subalgebra and $\fu$ is an ideal in $\fp$ 
follows easily from the relations 
\begin{equation}
  \label{eq:pararel}
(\Sigma^+ + \Sigma^+) \cap \cR \subeq \Sigma^+ 
\quad \mbox{ and } \quad 
(\Sigma^0 + \Sigma^+) \cap \cR \subeq \Sigma^+, 
\end{equation}
which are obvious consequence of \eqref{eq:pardef}. 
The subalgebra $\fu$ is closed and nilpotent. 
In fact, if $x_\Sigma$ is chosen as in Remark~\ref{rem:para-pos} 
and $N > \max \Sigma^+(x_\Sigma)$, then 
$N$-fold brackets in $\fu$ vanish. Likewise 
$$\fn= \sum_{\alpha \in -\Sigma^+} \g_\alpha $$
is a nilpotent closed subalgebra. As a Banach space, $\g$ has 
the direct sum decompositions  
\begin{equation}
  \label{eq:par-deco}
\g = \fn \oplus \fl \oplus \fu = \fn \oplus \fp. 
\end{equation}

\begin{rem} \label{rem:1.6}
Note that the subalgebra $\fl$ is weakly root graded with 
respect to the root system $\Delta^0 := 
\Sigma^0\cap \Delta$ and the grading subalgebra 
$$\g_{\Delta^0} := \h_{\Delta^0} 
+ \sum_{\alpha \in \Delta^0} \g_{\Delta, \alpha},  
\quad \mbox{ where } \quad 
\fh_{\Delta_0} := \Spann \{ \check \alpha \: \alpha \in \Delta^0\}. $$
The intersection $\fp_\Delta := \fp \cap \g_\Delta$ 
is a parabolic subalgebra of $\g_\Delta$ satisfying 
$$ \fp_\Delta = \fu_\Delta \rtimes \fl_\Delta \quad \mbox{ for } \quad 
\fl_\Delta := \fl \cap \g_\Delta \quad \mbox{ and } \quad 
\fu_\Delta := \fu \cap \g_\Delta. $$
We also put $\fn_\Delta := \g_\Delta \cap \fn$ and obtain the direct sum 
decomposition $\g_\Delta= \fn_\Delta \oplus \fp_\Delta$. 
\end{rem}

If $\g$ is a Lie algebra and $E \subeq \g$ a subspace, we write 
$$ \fn_\g(E) := \{ x \in \g \: [x,E] \subeq E \} $$
for the {\it normalizer of $E$ in $\g$}. 

\begin{lemma}\label{lem:normal} For any parabolic subalgebra 
$\fp = \fu \rtimes \fl$,  we have 
$$\fn_{\g}(\fp)=\fp \quad \mbox{ and } \quad 
\fn_\g(\fl) = \fl.$$
\end{lemma}

\begin{proof} (cf.\ \cite[Lemma~1.2.11]{Mue06})
Since $\fp$ and $\fl$ are $\h$-invariant, their normalizers 
contains $\h$, and hence are adapted to the root decomposition. 

Let $\alpha \in \cR \setminus \Sigma$ and 
$0 \not= y_\alpha \in \g_\alpha$.
{}From $y_\alpha \in [\fh,y_\alpha]$ it follows that 
$y_\alpha$ does not normalize $\fp$, which leads to $\fn_\g(\fp) = \fp$. 
We likewise derive that $\fn_\g(\fl) = \fl$. 
\end{proof}

\begin{remark}
In the case when $\cR=\Delta$ we can also derive the
equality $\fn_{\g}(\fu)=\fp$ from the representation theory of 
$\fsl_2(\C)$. This does however not hold in general. Consider 
for example the semidirect product $\fsl_2(\C) \ltimes \C^2$
given by the standard representation of $\fsl_2(\C)$. Here
$\cR=\{\pm \alpha, \pm2 \alpha \}$, and for the parabolic system
$\Sigma=\{\alpha, 2\alpha \}$ we have $\fu=\g_{\alpha} \oplus \g_{2\alpha}$.
However
$\fn_{\g}(\fu)=\g_{-\alpha} \oplus \h \oplus \g_{\alpha} \oplus \g_{2\alpha}$
since $[\g_{-\alpha}, \g_{\alpha}]={0}$.
\end{remark}

In the preceding lemma, we have seen that the closed Lie subalgebra 
$\fp$ of $\g$ is self-normalizing. From Lemmas~IV.11/12 in 
\cite{Ne04} we thus obtain that 
$$ N_G(\fp) := \{ g \in G \: \Ad(g)\fp = \fp\} $$
is a Lie subgroup of $G$ whose Lie algebra is 
$$ \L(N_G(\fp)) = \fn_\g(\fp) = \fp. $$

\begin{definition} \label{def:1.6} 
In the following we call any open subgroup $P$ of $N_G(\fp)$ 
a {\it (standard) parabolic subgroup} of $G$. Then 
$P_0 = \la \exp \fp \ra$
is a connected Lie subgroup of $G$, and since $\fn$ is a closed complement 
of $\fp$, we see that $P$ is a split Lie subgroup 
(cf.\ \cite[Def.~IV.6]{Ne04}). Therefore the quotient space 
$G/P$ carries a natural complex manifold structure for which 
the quotient map 
$q \: G \to G/P, g \mapsto gP$
is a submersion with holomorphic local cross section. 
\end{definition}

\begin{proposition} \label{prop:1.9} 
For any parabolic subgroup $P$ of $G$, we put 
$$ U := \exp_G \fu, \quad N  := \exp_G \fn \quad \mbox{ and } \quad 
L := P \cap N_G(\fl) \cap N_G(\fn). $$
Then the following assertions hold: 
\begin{description}
\item[\rm(a)] $U$ is a Lie subgroup with Lie algebra $\fu$ whose exponential 
function is a diffeomorphism. 
\item[\rm(b)] $N$ is a Lie subgroup with Lie algebra $\fn$ whose exponential 
function is a diffeomorphism. 
\item[\rm(c)] $L$ is a Lie subgroup with Lie algebra $\fl$ which 
acts holomorphically by conjugation on $N$ and $U$. 
\item[\rm(d)] The multiplication map $U \rtimes L \to P, (u,l) \mapsto 
ul$  
is biholomorphic onto an open subgroup of $P$.
\item[\rm(e)] The multiplication map $N \times L \times U \to G$ 
is biholomorphic onto an open subset of $G$.  
\end{description}
\end{proposition}

\begin{proof} (a) Let $*$ denote the Baker--Campbell--Hausdorff 
multiplication on the nilpotent Lie algebra $\fu$, turning it into a 
simply connected complex Banach--Lie group $(\fu,*)$. 
Then the exponential function 
$$ \exp \: \fu \to G, \quad x \mapsto \exp_G x $$
is the unique morphism of Banach--Lie groups $(\fu,*) \to G$ 
integrating the inclusion map $\fu \to \g$. This implies that 
$U = \exp_G(\fu)$ is a subgroup of $G$. 

To see that it is a Lie subgroup, we recall that $P$ is a Lie subgroup 
of $G$ (Definition~\ref{def:1.6}). 
This implies that 
$$ U_1 := \{ p \in P \: (\Ad(p)-\1)(\fp) \subeq \fu\} $$
also is a Lie subgroup of $P$ and hence of $G$ 
(use Lemma~IV.11 in \cite{Ne04} and generalize 
Lemma~IV.12 slightly). Clearly, 
$$\L(U_1) = \{ x \in \fp \: [x,\fp] \subeq \fu\} \supeq \fu, $$
but, conversely, $[x,\fh] \subeq \fu$ implies $x \in \fu$, so that 
$\L(U_1) = \fu$. We conclude that  $U$ is the identity component of the 
Lie subgroup $U_1$, hence a Lie subgroup. 

Finally, we show that $\exp_U$ is a diffeomorphism. Since $U$ is connected 
and nilpotent, it suffices to show that $\exp_U$ is injective. 
So let $x \in \fu$ with $\exp_G x = \1$. Then 
$\1 = \Ad(\exp_G x) = e^{\ad x}$, and since $\ad x$ is nilpotent, 
we get $\ad x= 0$. Since the root decomposition yields 
$\z(\g) \subeq \g_0$, we see that $x = 0$. 

(b) follows from (a), applied to the parabolic subalgebra defined by 
$-\Sigma$. 

(c) That $L$ is a Lie subgroup follows from \cite[Lemmas~IV.11/12]{Ne04}. 
In view of Lemma~\ref{lem:normal}, its Lie algebra is 
$\fp \cap \fn_\g(\fl) \cap \fn_\g(\fn) = \fl$. 

Since $L$ normalizes $\fu$ and $\fn$, its acts by the adjoint action 
on these Lie algebras, and (a) and (b) imply that this corresponds to a 
holomorphic action by conjugation on the groups $U$ and $N$. 

(d) Since $L$ acts holomorphically by conjugation on $U$, the 
semidirect product $U \rtimes L$ is a complex Banach--Lie group. 
As $\fp = \fu \rtimes \fl$ is a semidirect sum, 
the canonical homomorphism 
$U \rtimes L \to P, (u,l) \mapsto ul$ 
is a morphism of Banach--Lie groups 
which is a local diffeomorphism. Therefore $UL$ is an open subgroup 
of $P$ and the multiplication map $U \times L \to UL$ 
is a covering. It remains to see that this map is injective, i.e., 
that $U \cap L = \{\1\}$. 

Let $x \in \fu$. If $\exp_G x \in L$, then the unipotent operator 
$e^{\ad x}$ normalizes $\fl$, and this implies that 
$\ad x = \log(e^{\ad x})$ also preserves $\fl$. Hence 
$x \in \fn_\g(\fl) = \fl$ leads to $x =0$. 
 
(e) The direct product Lie group 
$N \times (U \rtimes L)$ acts smoothly on $G$ by 
$$ (n,(u,l)).g := n g l^{-1} u^{-1}, $$
so that (e) means that the orbit map of $\1$ is a diffeomorphism 
onto an open subset. That the differential in $\1$ is a linear 
isomorphism follows from equation~\eqref{eq:par-deco}, so that the 
orbit map is a local diffeomorphism by the Inverse Function Theorem, 
hence a covering map whose range is an open subset of~$G$. 
It remains to verify that the orbit map is injective, i.e., 
$P \cap N = \{\1\}$. 
This can be seen as in (d): If the $\ad$-nilpotent element 
$x \in \fn$ satisfies 
$e^{\ad x} \fp = \fp$, then $x \in \fn_\g(\fp) = \fp$ yields $x =0$, 
and this prove that $N \cap P=\{\1\}$. 
\end{proof}

\section{Coinduced Banach representations} \label{sec:2}

In this section we discuss coinduced representation of Banach--Lie 
algebras. After introducing a suitable 
notion of continuity for linear maps $\alpha \: U(\g) \to E$ 
on the enveloping algebra $U(\g)$ of a Banach--Lie algebra $\g$ with values 
in a Banach space $E$, we show that the space 
$\Hom(U(\g),E)_c$ of all these maps carries a natural Fr\'echet 
space structure for which the left and right action of $\g$ 
are continuous. If $\fp \subeq \g$ is a complemented closed subalgebra, 
this leads for each continuous representation $\rho \: \fp \to \gl(E)$ 
to a Fr\'echet structure on the corresponding coinduced representation 
$\Hom_\fp(U(\g),E)_c$. Specializing all that to a parabolic subalgebra 
$\fp$ of a weakly root graded Lie algebra $\g$, our main result 
is that, provided $E$ decomposes into finitely many $\fh$-weight 
spaces, the subspace  $\Hom_\fp(U(\g),E)_c^{[\g_\Delta]}$ of 
$\g_\Delta$-finite elements 
is closed and a Banach space (Theorem~\ref{thm:weightfin-alg}). 
As we shall see in the next section, this implies that on the group 
level, holomorphic induction from Banach representations of $P$ 
yields Banach representations of $G$.

\subsection*{Coinduced Fr\'echet representations} 

\begin{defn} Let $V$ and $W$ be Banach spaces. 
For an $n$-linear map $m \: V^n \to W$, we define 
$$ \|m\| := \sup \{ \|m(v_1,\ldots, v_n)\| \: v_1, \ldots, v_n \in V, 
\|v_i\| \leq 1\} 
\in [0,\infty]. $$
Then $m$ is continuous if and only if $\|m\| < \infty$, and we have 
$$ \|m(v_1,\ldots, v_n)\| 
\leq \|m\| \|v_1 \| \cdots \|v_n\|, \quad v_i \in V. $$
We thus obtain on the space $\Mult^n(V,W)$ of continuous 
$W$-valued $n$-linear 
maps on $V^n$ the structure of a Banach space. 
\end{defn}

\begin{defn}
Let $\g$ be a Banach--Lie algebra and $E$ be a Banach space. 
We call a linear map 
$\beta \: U(\g) \to E$ {\it continuous} if for each $n \in \N_0$,  
the $n$-linear map 
$$ \beta_n \: \g^n \to E, \quad (x_1,\ldots, x_n) \mapsto 
\beta(x_1\cdots x_n) $$
is continuous.  We write $\Hom(U(\g),E)_c$ for the set of all 
continuous linear maps and define a family of seminorms 
$$ p_n \: \Hom(U(\g),E)_c \to \R,\quad 
\beta \mapsto \|\beta_n\|. $$
\end{defn} 

\begin{lem} \label{lem:frechet} 
With respect to the sequence $(p_n)_{n \in \N_0}$ of 
seminorms, the space $\Hom(U(\g),E)_c$ is a Fr\'echet space 
for which all evaluation maps 
$$ \ev_D \: \Hom(U(\g),E)_c \to E, \quad \beta \mapsto \beta(D) $$
are continuous. 
Moreover, the left representation of $U(\g)$ on this space 
defined by 
$$ \g \times \Hom(U(\g),E)_c \to \Hom(U(\g),E)_c, \quad 
(x,\beta) \mapsto -\beta \circ \lambda_x, \quad 
\lambda_x(D) = xD, $$
and the right representation 
$$ \g \times \Hom(U(\g),E)_c \to \Hom(U(\g),E)_c, \quad 
(x,\beta) \mapsto \beta \circ \rho_x, \quad 
\rho_x(D) = Dx, $$
are continuous bilinear maps. 
\end{lem} 

\begin{prf} Since the canonical map 
$\Hom(U(\g),E)_c \to \prod_{n \in \N_0} \Mult^n(\g,E)$
is injective, the sequence of seminorms $(p_n)_{n \in \N_0}$ defines 
on $\Hom(U(\g),E)_c$ a metrizable locally convex topology. It also follows 
immediately from the definitions, that all evaluation maps $\ev_D$ 
are continuous because we may write 
$D$ as a finite sum 
$D = D_0 + \ldots + D_k$, where 
$D_j$ is a $j$-fold product $x_1\cdots x_j$ of elements $x_i \in \g$. 

To verify that $\Hom(U(\g),E)_c$ is complete, i.e., a Fr\'echet space, let 
$(\beta^{(n)})_{n \in \N}$ be a Cauchy sequence in $\Hom(U(\g),E)_c$. 
Then, for each $k \in \N_0$, the sequence $(\beta^{(n)}_k)_{n \in \N}$ 
define a Cauchy sequence in the Banach space $\Mult^k(\g,E)$, which 
converges to some limit $\beta_k$. We also note that for each $D \in U(\g)$, 
the sequence $(\beta^{(n)}(D))_{n \in \N}$ in $E$ is Cauchy, hence 
convergent. If we write $\beta \: U(\g) \to E$ for the pointwise 
limit, we thus obtain a linear map whose corresponding $k$-linear maps 
$\g^k \to E$ coincide with the continuous maps $\beta_k$. This implies 
that $\beta \in \Hom(U(\g),E)_c$ and that $\beta^{(n)} \to \beta$ 
holds in the metric topology on $\Hom(U(\g),E)_c$. 

Next we observe that for $\beta \in \Hom(U(\g),E)_c$ and 
$x \in \g$ we have 
\begin{equation}
  \label{eq:semiest}
p_n(\beta \circ \lambda_x) \leq \|x\| \cdot p_{n+1}(\beta) 
\quad \mbox{ and } \quad 
p_n(\beta \circ \rho_x) \leq \|x\|\cdot p_{n+1}(\beta). 
\end{equation}
This proves that the corresponding bilinear maps 
$$ \g \times \Hom(U(\g),E)_c \to \Mult^n(\g,E), \quad 
(x,\beta) \mapsto (\beta \circ \lambda_x)_n, (\beta \circ \rho_x)_n, $$
are continuous. Since the topology on $\Hom(U(\g),E)_c$ is obtained by 
the embedding into the product of the spaces $\Mult^n(\g,E)$, 
the assertion follows. 
\end{prf}

Next we assume that $\fp$ and $\fn$ are closed subalgebras of $\g$ with 
$\g = \fn \oplus \fp$ (topological direct sum) and that 
$(\rho_E, E)$ is a continuous representation of $\fp$ on the 
Banach space $E$. 

\begin{defn} For a continuous morphism $\rho_E \: \fp \to \gl(E)$ 
of Banach--Lie algebras, we consider the subspace 
$$ \Hom_\fp(U(\g), E)_c 
:= \{ \beta \in \Hom(U(\g),E)_c \: 
(\forall x \in \fp)\, \beta \circ \rho_x = - \rho_E(x) \circ \beta\}. $$
{}From Lemma~\ref{lem:frechet} we immediately derive that 
this is a closed subspace of the Fr\'echet space 
$\Hom(U(\g),E)_c$, hence also a Fr\'echet space, and that the left 
action of $\g$ on this space is continuous. 
It is called the {\it coinduced representation} defined by $(\rho_E,E)$. 
\end{defn}

Since the multiplication map 
$U(\fn) \otimes U(\fp) \to U(\g)$ is a linear isomorphism 
by the Poincar\'e--Birkhoff--Witt Theorem, 
the restriction map 
$$ R \: \Hom_\fp(U(\g), E)_c \cong \Hom(U(\fn), E)_c $$
is injective. We even have more: 

\begin{lem} \label{lem:2.5} 
$R$ is a topological isomorphism of Fr\'echet spaces. 
\end{lem}

\begin{prf} Clearly, $R$ is injective, and the PBW Theorem implies that 
it is bijective on the algebraic level. Each 
$\alpha \in \Hom(U(\fn),E)$ extends by 
$$ \tilde \alpha(D \cdot x_1 \cdots x_n) := 
(-1)^n \rho_E(x_n) \cdots \rho_E(x_1) \alpha(D), \quad 
x_i \in \fp, D \in U(\fn), $$
to a unique linear map $\tilde\alpha \in \Hom_\fp(U(\fg),E)$. 
It remains to show that the associated $n$-linear maps $\tilde\alpha_n$ 
are continuous if all maps $\alpha_n$ are continuous. 

To this end, we show that there exist constants 
$C_{j,n}$, $j \leq n$, such that 
\begin{equation}
  \label{eq:esti}
 p_n(\tilde \alpha) \leq \sum_{j \leq n} C_{j,n} p_j(\alpha). 
\end{equation}
To verify these estimates, let $C_\g\geq 0$ be such that 
$$ \|[x,y]\| \leq C_\g \|x\| \cdot \|y\| \quad \mbox { for } 
\quad x,y \in \g, $$
and assume w.l.o.g.\ that 
$$ \| x + y\| = \max(\|x\|,\|y\|) \quad \mbox {for } \quad x \in \fp, y \in \fn.$$
We now prove \eqref{eq:esti} by induction on $n$. For $n = 0$, 
we only deal with the constant term, so that 
$p_0(\tilde\alpha) = p_0(\alpha)$. Next we assume that 
$n > 0$ and let $z_i = x_i + y_i \in \g$, $x_i \in \fp$, $y_i \in \fn$, 
$i =1,\ldots, n$ with $\|z_i \| \leq 1$. Then 
\begin{align*}
&\tilde\alpha(z_1 \cdots z_n) 
= \tilde\alpha(z_1 \cdots z_{n-1} \cdot (x_n + y_n)) \\
&=  - \rho_E(x_n)\tilde\alpha(z_1 \cdots z_{n-1}) 
-  \tilde\alpha([y_n, z_1 \cdots z_{n-1}]) 
+  \tilde\alpha(y_n \cdot z_1 \cdots z_{n-1}) \\
&=  - \rho_E(x_n)\tilde\alpha(z_1 \cdots z_{n-1}) 
-  \tilde\alpha([y_n, z_1 \cdots z_{n-1}]) 
+  (\alpha \circ \lambda_{y_n})\,\tilde{}\,(z_1 \cdots z_{n-1}). 
\end{align*}
For the last term we obtain with the induction hypothesis and 
\eqref{eq:semiest} the estimate 
$$ \|(\alpha \circ \lambda_{y_n})\,\tilde{}(z_1 \cdots z_{n-1})\| 
\leq p_{n-1}\big((\alpha \circ \lambda_{y_n})\,\tilde{}\big) 
%\leq \sum_{j \leq n-1} C_{j,n-1} p_j(\alpha \circ \lambda_{y_n}) 
\leq \sum_{j \leq n-1} C_{j,n-1} p_{j+1}(\alpha). $$
 For the other two terms, we have 
 \begin{align*}
\|\rho_E(x_n)\tilde\alpha(z_1 \cdots z_{n-1}) \| 
&\leq p_{n-1}(\rho_E(x_n)\circ \tilde \alpha) 
\leq \|\rho_E\| p_{n-1}(\tilde \alpha) \\
&\leq \|\rho_E\| \sum_{j \leq n-1} C_{j,n-1} p_j(\alpha)  
 \end{align*}
and 
$[y_n, z_1 \cdots z_{n-1}] 
= \sum_{i = 1}^{n-1} z_1 \cdots z_{i-1} [y_n, z_i] z_{i+1} \cdots z_{n-1}$
leads to 
$$ \|\tilde\alpha([y_n, z_1 \cdots z_{n-1}])\| 
\leq (n-1) C_\g p_{n-1}(\tilde\alpha)
\leq (n-1) C_\g \sum_{j \leq n-1} C_{j,n-1} p_j(\alpha). $$
This completes the inductive proof of \eqref{eq:esti}. 
\end{prf}

\begin{rem} \label{rem:2.6} Let $k \in \N_0$ and consider in 
$\Hom(U(\fn),E)_c$ the subspace 
$$\Hom(U(\fn),E)_c^k := \{ \alpha \in \Hom(U(\fn),E)_c \: 
(\forall n \geq k)\, \alpha_n = 0\}. $$
This is a closed subspace of the Fr\'echet space 
$\Hom(U(\fn),E)_c$, and since all but finitely many of the 
seminorms $p_n$ vanish on this space, it actually is a Banach space. 
\end{rem}

\subsection*{Applications to weakly root graded Lie algebras} 

In this subsection we return to the setting where 
$\fp$ is a parabolic subalgebra of the weakly root graded 
complex Banach--Lie algebra $\g$.

\begin{lemma} \label{lem:fin-iso} 
For a $\g$-module $V$ with an $\h$-weight decomposition,  
the following are equivalent: 
\begin{description}
\item[\rm(a)] The set of $\h$-weights of $V$ is finite. 
\item[\rm(b)] $V$ is a locally finite $\g_\Delta$-module with 
finitely many isotypic components.  
\end{description}
Then $V$ is a semisimple $\g_\Delta$-module, hence in particular 
a direct sum of simple ones. 
\end{lemma}

\begin{proof} (a) implies that the $\g_\Delta$-module $V$ is integrable 
in the sense that each root vector $x_\alpha \in \g_{\Delta,\alpha}$ acts as a 
nilpotent operator on $V$. Hence \cite[Thm.~A.1(1), Prop.~A.2]{Ne03} 
imply that $V$ is a locally finite (and therefore semisimple) 
$\g_\Delta$-module 
with finitely many isotypic components. That (b) implies (a) is trivial. 
\end{proof}

\begin{lemma} \label{lem:2.8} Let $\fk \subeq \g$ be a subalgebra 
for which $\g$ is a locally finite $\fk$-module. 
Then for any representation $(\pi,V)$ of $\g$, 
the subspace 
$$V^{[\fk]} := \{ v \in V \: \dim \big(U(\fk)v\big) < \infty\} $$
of $\fk$-finite vectors is a $\g$-submodule. 
\end{lemma} 

\begin{proof} 
%For any $\Delta$-graded Lie algebra $\g$, the space 
%$\g$ is a $\g_\Delta$-finite representation under the adjoint action 
%(Lemma~\ref{lem:fin-iso}). 
Since $\g$ is locally $\fk$-finite, the tensor algebra over 
$\g$, and hence $U(\g)$, is a locally $\fk$-finite $\g$-module. 
This property carries over to the tensor product 
$U(\g) \otimes V^{[\fk]}$, and hence to its image under the 
$\fk$-equivariant evaluation map $U(\g) \otimes V^{[\fk]} \to V$. 
We conclude that $U(\g)V^{[\fk]} \subeq V^{[\fk]}$. 
\end{proof}

\begin{defn} We define the {\it Weyl group} 
${\cal W} := \cW(\g_\Delta, \fh) \subeq \GL(\fh^*)$ as the subgroup 
generated by 
the reflections $r_\alpha(\beta) := \beta - \beta(\check \alpha) \alpha$. 
\end{defn}

\begin{thm}{\rm(Finiteness Theorem)} \label{thm:weightfin-alg} 
Let $\fp \subeq \g$ be a parabolic subalgebra and 
$(\rho,E)$ be a $\fp$-module decomposing into a 
finite sum of $\h$-weight spaces 
and $V := \Hom_\fp(U(\g),E)^{[\fh]}$ be the $\fh$-finite part of 
the corresponding coinduced representation. 
Then the following assertions hold: 
\begin{description}
\item[\rm(a)] $V$ is a direct sum of $\fh$-weight spaces   
and the set $\cP(V,\fh)$ of $\fh$-weights in $V$ is contained in 
$\cP(E,\fh)- \Spann_{\N_0} \Sigma^-$. 
\item[\rm(b)]  The $\g_\Delta$-finite subspace 
$V^{[\g_\Delta]}$ is a $\g$-module decomposing 
into finitely many $\h$-weight spaces. 
\item[\rm(c)] $V$ is a locally finite module of 
$\fn_\Delta \rtimes \fl_\Delta$. 
\item[\rm(d)] There exists a $k \in \N$ with 
$V^{[\g_\Delta]} = \{ v \in V \: \fu_\Delta^k.v = \{0\}\}.$
\item[\rm(e)] There exists a $k \in \N$ with 
$V^{[\g_\Delta]} \subeq \{ v \in V \: \fn^k.v = \{0\}\}. $
\end{description}
\end{thm}

\begin{proof} (a) First we recall that the multiplication map 
$U(\fn) \otimes U(\fp) \to U(\g)$ is a linear isomorphism 
(Poincar\'e--Birkhoff--Witt), 
so that we obtain a linear isomorphism 
$\Hom_\fp(U(\g), E) \cong \Hom(U(\fn), E),$
of $\h$-modules. If $E = \bigoplus_{\alpha \in \cP(E,\fh)} E_\alpha$ 
is the finite weight space decomposition of $E$, we accordingly obtain a 
product decomposition of $\fh$-modules 
$$ \Hom(U(\fn), E) 
\cong \prod_{\alpha \in \cP(E,\fh)} \Hom(U(\fn), E_\alpha). $$
Since $\fn$ is a finite sum of $\fh$-weight spaces, the enveloping algebra 
$U(\fn)$ is an (infinite) direct sum of $\fh$-weight spaces 
$U(\fn)_\beta$, which leads to a product decomposition of 
$\fh$-modules 
$$ \Hom(U(\fn), E_\alpha) \cong \prod_\beta \Hom(U(\fn)_\beta, E_\alpha), $$
and 
$\Hom(U(\fn)_\beta, E_\alpha) \subeq V_{\alpha - \beta}.$
Therefore the set $\cP(V,\beta)$ of $\fh$-weights on $V$ is given by 
$$ \cP(V,\fh) = \cP(E,\fh) - \Spann_{\N_0} \Sigma^-, $$ 
and $V = \sum_\alpha \sum_\beta \Hom(U(\fn)_\beta, E_\alpha)$ follows from the 
fact that an element of 
$\Hom(U(\fn), E_\alpha)$ is $\fh$-finite if and only if only finitely 
many components in the spaces $\Hom(U(\fn)_\beta, E_\alpha)$ are non-zero. 

(b) Since $\g$ is a locally finite module of the subalgebras $\fh$ and 
$\g_\Delta$, Lemma~\ref{lem:2.8} implies that $V$ is a $\g$-submodule  
of $\Hom_\fp(U(\g),E)$ and that $V^{[\g_\Delta]}$ is a $\g$-submodule of 
 $V$. 

Let $\cR^+$ be a positive system contained in $-\Sigma$ 
(Definition~\ref{def:parasys}), so 
that we necessarily have $\Sigma^- \subeq \cR^+$. We have seen in (a) 
that 
\begin{equation}
  \label{eq:firstinc}
 \cP(V,\fh) \subeq \cP(E,\fh) - \Spann_{\N_0}\Sigma^- 
\subeq \cP(E,\fh) - C, 
\end{equation}
where $C := \Spann_{\R_+} \cR^+$ is a polyhedral convex cone 
containing no non-zero linear subspace. 

The local finiteness of the $\g_\Delta$-module $V^{[\g_\Delta]}$ 
implies that the 
set $\cP_f := {\cal P}(V^{[\g_\Delta]},\fh)$ is invariant under the Weyl group 
${\cal W}$ (\cite{MP95}), which leads to 
\begin{equation}
  \label{eq:weightinc}
 \cP(V^{[\g_\Delta]},\fh) \subeq D := \bigcap_{w \in \cW} w(\mbox{conv}(\cP(E,\fh)) - C).
\end{equation}
The set $D$ is convex and we claim that it is compact. 
If this is not the case, then there exists a non-zero 
$\beta \in \fh_\R^*$ with $D + \R_+ \beta \subeq D$ 
(\cite[Prop.~V.1.6]{Ne00}). Pick $d \in D$. Then 
$d + \R_+ \beta \subeq w(\mbox{conv}(\cP(E,\fh)) - C)$ implies that 
$\beta \in - wC$ for each $w \in \cW$ (\cite[Prop.~V.1.6]{Ne00}). 
We thus arrive at $\cW \beta \subeq - C$. 
By definition of a positive system, there exists an element 
$x \in \fh$ with $\cR^+(x) > 0$, and then 
$\gamma(x) > 0$ holds for each non-zero element $\gamma \in C$, 
in particular $(w\beta)(x) < 0$ for each $w \in \cW$ and for 
$\gamma := \sum_{w \in \cW} w\beta$ we also obtain 
$\gamma(x) < 0$. Since $\gamma$ is $\cW$-invariant, it 
vanishes on all coroots, which leads to the contradiction $\gamma = 0$. 
We conclude that 
\begin{equation}
  \label{eq:intersect}
\bigcap_{w \in W} w C = \{0\}.  
\end{equation}

This proves that $D$ is compact. Now the discreteness of the 
set \break ${\cal P}(E,\fh) + \Spann_\Z\cR$ in $\fh_\R^*$ 
implies that ${\cal P}_f$ is finite. As $\fh$ is diagonalizable on $V$, 
(b) follows. 

(c) In view of the PBW Theorem, it suffices to show that 
$V$ is a locally finite module of $\fn_\Delta$ and $\fl_\Delta$. 

To see that it is locally finite for $\fn_\Delta$, we pick 
$x \in \fh$ with $\Sigma^0(x) = \{0\}$ and $\Sigma^-(x) \geq 1$ 
(cf.\ Remark~\ref{rem:para-pos}). Then (a) implies that 
$\cP(V,\fh)(x)$ is bounded from above. 
Therefore $\g_\alpha V_\beta \subeq V_{\alpha + \beta}$ 
and $\Sigma^-(x) \geq 1$ imply that for each weight $\beta$ of $V$, 
there exists a $k \in \N$ with 
$V_{\beta + \alpha_1 + \ldots + \alpha_k} = \{0\}$
for $\alpha_i \in \Sigma^-$. We thus derive that $V$ is a locally 
nilpotent $\fn$-module, hence in particular locally finite for the 
finite dimensional Lie algebra $\fn_\Delta$. 

Next we note that on each space 
$\Hom(U(\fn)_\beta, E_\alpha)$, the eigenvalue of 
$x \in \fh$ is given by $\alpha(x) - \beta(x)$ and, for each $c \in \R$, 
the set 
$\{ \beta \in \Spann_{\N_0} \Sigma^- \: \beta(x) = c \}$
is finite. Since $x$ commutes with $\fl$, the eigenspaces of $x$ in 
$V$ are $\fl$-invariant, and the preceding argument shows that they  
decompose into finitely many $\fh$-weight spaces. 
Therefore Lemma~\ref{lem:fin-iso} implies that 
each such eigenspace is a locally finite $\fl_\Delta$-module. 
Here we use Remark~\ref{rem:1.6} to see that 
$\fl$ is $\Delta^0$-graded, so that Lemma~\ref{lem:fin-iso} 
applies. 

(d) In view of (c) and 
$U(\g_\Delta) = U(\fu_\Delta) U(\fl_\Delta)U(\fn_\Delta)$, 
the subspace $V^{[\g_\Delta]}$ of $\g_\Delta$-finite elements 
coincides with the subspace $V^{[\fu_\Delta]}$ of $\fu_\Delta$-finite 
elements. 
Further, the finiteness of $\cP_f$ and 
$\Sigma^+(x) \leq -1$ entail the existence of some $k \in \N$ 
with 
$$ V^{[\g_\Delta]} \subeq \{  v \in V \: (\fu_\Delta)^k.v = \{0\}\}. $$
Conversely, the condition $(\fu_\Delta)^k.v = \{0\}$ obviously 
implies that $v$ is $\fu_\Delta$-finite and hence $\g_\Delta$-finite 
by (c). 
We therefore obtain the desired equality (d). 

(e) follows from the finiteness of $\cP_f$ and 
$\Sigma^-(x) \geq 1$. 
\end{proof}

\begin{rem} \label{rem:triv-mod} 
Equations \eqref{eq:weightinc} and \eqref{eq:intersect} 
have an interesting 
consequence. If $E$ is a trivial $\fp$-module, then $\cP(E,\fh) = \{0\}$ 
implies that 
$$\cP(V^{[\g_\Delta]},\fh)
 \subeq - \bigcap_{w \in \cW} wC = \{0\},$$ 
so that $\fh$ acts trivially on $V^{[\g_\Delta]}$. This in turn 
implies that all root spaces $\g_\alpha$, $\alpha \in \Delta$, 
act trivially on $V^{[\g_\Delta]}$. 
If $\g$ is root graded, we conclude that the $\g$-module 
$V^{[\g_\Delta]}$ is trivial, and we find that 
$$ V^{[\g_\Delta]} = \Hom_\fp(U(\g),E)^\g \cong E. $$
\end{rem}

The preceding theorem is of an algebraic nature, but it 
has some interesting consequences for Banach representations. 

\begin{thm}{\rm(Coinduced Banach representations)} \label{thm:coind-ban}
Let $(\rho_E,E)$ be a Banach representation of $\fp$ 
which decomposes into finitely many $\fh$-weight spaces 
and endow $V_c := \Hom_\fp(U(\g),E)_c$ with its natural Fr\'echet structure. 
Then the subspace $V_f := V_c^{[\g_\Delta]}$ of $\g_\Delta$-finite elements 
in $V_c$ is a closed $\g$-invariant subspace which is Banach 
and the representation of $\g$ on this space is continuous. 
\end{thm} 

\begin{prf} First we apply Theorem~\ref{thm:weightfin-alg} 
to the subspace $V_c$ of $\Hom_\fp(U(\g),E)$. Since each 
$\g_\Delta$-finite vector is in particular $\fh$-finite, 
$V_f \subeq \Hom_\fp(U(\g),E)^{[\fh]}$, so that 
Theorem~\ref{thm:weightfin-alg}(d) implies the existence of a $k \in \N$ 
with 
$$ V_f = \{ v \in V_c \: (\fu_\Delta^k).v = \{0\}\}. $$
This implies that $V_f$ is a closed subspace of the Fr\'echet space 
$V_c$. Next we use Theorem~\ref{thm:weightfin-alg}(e) to find a 
$k \in \N$ with $\fn^k.V_f = \{0\}$. 
Then, for each $\alpha \in V_f$ and $m \in \N$, the relation 
$$ (\fn^k.\alpha)(\fn^m) = \alpha(\fn^{k+m}) = \{0\} $$
implies that 
$$ V_f\res_{U(\fn)} \subeq \Hom(U(\fn),E)_c^k, $$
and since the right hand side is a Banach space (Remark~\ref{rem:2.6}),  
Lemma~\ref{lem:2.5} implies that $V_f$ is a Banach space. 
The continuity of the $\g$-action on the Banach space 
$V_f$ is inherited from the continuity on the Fr\'echet space $V_c$. 
\end{prf}

We now turn to some additional information on the 
$\g$-representation on~$V_f$. 

\begin{prop} \label{prop:irred} 
Let $(\rho_E,E)$ be an {\em irreducible} Banach representation of $\fp$ 
which decomposes into finitely many $\fh$-weight spaces and assume that 
$V_f$ is non-zero. Then the following assertions hold: 
\begin{description}
\item[\rm(i)]  $\fu$ acts trivially on $E$. 
\item[\rm(ii)]  Each closed 
$\g$-submodule of the Banach space $V_f$ contains the space 
$E \cong \Hom_\fp(U(\g),E)_c^\fn$ of $\fn$-invariants. In particular, 
$E$ generates a unique closed minimal submodule which is irreducible.  
\item[\rm(iii)]  For each weight $\alpha \in \cP(E,\fh)$, we have 
$$V_{f,\alpha} = E_\alpha \quad \mbox{ and } \quad 
\alpha(\check \beta) \in \N_0 \quad \mbox{ for } \quad 
\beta \in \Sigma^- \cap \Delta. $$
If, in addition, $\dim E = 1$ and  
$\lambda := \rho\res_\fh$,  then the weight space  
$V_{f,\lambda}$ is one dimensional. 
\item[\rm(iv)] For each $\fp$-equivariant continuous map 
$\phi \: V_f \to E$, there exists a unique $\psi \in \End_\fp(E)$ with 
$\phi = \psi \circ \ev_\1.$
\end{description}
\end{prop}

\begin{prf} (i) Since $E$ decomposes into finitely many $\fh$-weight spaces, 
we have $\fu^k.E = \{0\}$ for some $k$. In particular, 
$E^\fu$ is a non-zero closed $\fp$-submodule. Since $E$ was assumed to be 
irreducible, it follows that $\fu \subeq \ker \rho_E$, and hence that 
$E$ is an irreducible $\fl$-module. 

(ii) We have already seen in Theorem~\ref{thm:weightfin-alg}(e)  
that $V_f$ is a nilpotent $\fn$-module, and this property is inherited 
by any $\g$-submodule $W \subeq V_f$. Hence 
$W^\fn \not=\{0\}$ whenever $W \not=\{0\}$. 
For the  $\fn$-invariant vectors we observe that  
$$\Hom_\fp(U(\g), E)_c^\fn \cong \Hom(U(\fn),E)_c^\fn 
\cong \Hom(U(\fn)/\fn U(\fn),E)_c \cong E. $$
In this sense we identify $E$ with the $\fl$-submodule 
$\Hom_\fp(U(\g),E)_c^\fn$ of \break $\Hom_\fp(U(\g),E)_c$. 
As $W^\fn$ is a non-zero closed $\fl$-submodule of $E$, we thus obtain 
$W^\fn = E$ and therefore $E \subeq W$. 

(iii) Next we pick $x_\Sigma \in \fh$ as in Remark~\ref{rem:para-pos}. 
Then $x_\Sigma$ is central in $\fl$ and $\rho_E(x_\Sigma)$ is diagonalizable. 
Since all eigenspaces of this element are $\fl$-submodules, the 
irreducibility of $E$ entails that $\rho_E(x_\sigma) = c \id_E$ for some 
$c \in \C$. On the other hand, $\Sigma^-(x_\Sigma) \leq -1$, and 
$U(\fn)$ decomposes into eigenspaces 
$U(\fn)_\mu$ of $x_\Sigma$ with $\mu \in ]-\infty, -1]$, where 
$U(\fn)_0 = \C \1$. Therefore 
$$ V_f \subeq \bigoplus_\mu \Hom(U(\fn)_\mu ,E), $$
and the eigenvalue of $x_\Sigma$ on  
$\Hom(U(\fn)_\mu ,E)$ is $c - \mu$, which is different from 
$c$ if $\mu\not=0$. This proves that 
$$V_{f,c}(x_\Sigma) = \{ v \in V_f \: x_\Sigma v = c v \} =  E,$$
and therefore 
$V_{f,\alpha} = E_\alpha$ for each $\fh$-weight $\alpha$ of $E$. 

Let $\beta \in \Delta \cap \Sigma^-$ 
and $\g_\Delta(\beta) := \g_{\Delta,\beta} + \g_{\Delta,-\beta} 
+ \C \check \beta \cong \fsl_2(\C)$ be 
the $\fsl_2$-subalgebra corresponding to 
$\beta$. Then the preceding argument shows that each element of 
$E_\alpha$ generates a finite dimensional $\g_\Delta(\beta)$-module 
for which the $\check\beta$-eigenvalues are contained in 
$\alpha(\check\beta) - 2 \N_0$, so that $\fsl_2$-theory yields 
$\alpha(\check \beta) \in \N_0$. 

(iv) Since $\fh \subeq \fp$, $\phi$ annihilates 
each weight space $V_{f,\beta}$ with $\beta \not\in \cP(E,\fh)$. 
In view of (iii), the fact that $V_f$ is a direct sum of $\fh$-weight spaces 
implies that 
$$ V_f \cong E \oplus \sum_{\beta\not\in \cP(E,\fh)} V_{f,\beta}, $$
and by identifying 
$\Hom_\fp(U(\g),E)_c$ with 
$\Hom(U(\fn),E)_c$ (Lemma~\ref{lem:2.5}), 
we see that $\ker(\ev_\1)$ coincides with the 
sum of all $\fh$-weight spaces $V_{f,\beta}$, $\beta \not\in \cP(E,\fh)$. 
We conclude that $\phi$ vanishes on $\ker(\ev_\1)$, hence induces a 
unique continuous $\fp$-equivariant map $\psi \in \End_\fp(E)$ with 
$\psi \circ \ev_\1 = \phi$. 
\end{prf}

\section{Holomorphically induced representations} \label{sec:3}

Let $G$ be a weakly root graded complex Banach--Lie group. 
Each closed subalgebra of $\g := \L(G)$ integrates to an 
integral subgroup (cf.\ \cite{GN09}). 
In particular, we obtain an integral subgroup $G_\Delta$ corresponding 
to the inclusion $\g_\Delta \into \g$. We write $H$ for the Cartan 
subgroup of $G_\Delta$ corresponding to the Cartan subalgebra $\h$ 
and $T \subeq H$ for its (unique) maximal compact subgroup, a real torus. 
For $r := \dim T$, we have $T \cong \T^r$, $H \cong (\C^\times)^r$,  
and $H \cong T_\C$ is the universal complexification of~$T$, so that 
we may identify the group $\hat T = \Hom(T,\T)$ of 
continuous characters of $T$ with the group 
$\hat H := \Hom_{\cal O}(H,\C^\times)$ of holomorphic characters of~$H$. 

\begin{theorem}{\rm(Peter--Weyl)} \label{thm:peter-weyl}
If $\pi \: K \to \GL(V)$ is a homomorphism 
defining a continuous linear action 
of the compact group $K$ on the complete locally convex space $V$, 
then the following assertions hold: 
\begin{description}
\item[\rm(1)] The space $V^{[K]}$ of $K$-finite vectors is dense in $V$. 
\item[\rm(2)] For each irreducible representation 
$[\chi] \in \hat K$, there is a continuous projection 
$p_{[\chi]} \: V \to V_{[\chi]}$ onto the corresponding isotypic subspace. 
\item[\rm(3)] If the support 
$\supp(\pi,V) := \{ [\chi] \in \hat K \: 
V_{[\chi]} \not=\{0\}\}$ is finite, then 
$V = \bigoplus_{[\chi] \in \hat K} V_{[\chi]}$
is a finite sum. 
\end{description}
\end{theorem}

\begin{proof} The first two assertions follows from the 
Peter--Weyl Theorem (\cite[Thm.~3.5.1]{HoMo98}). 
To derive (3), we note that 
$p := \sum_{[\chi] \in \supp(\pi,V)} p_{[\chi]}$
is a continuous projection whose range is $V^{[K]}$, so that (1) implies that 
$p$ is surjective. 
\end{proof} 

\begin{lemma} \label{lem:2.1} For each holomorphic Banach representation 
$(\pi,V)$ of the complex torus $H$, the set ${\cal P}(V,\fh)$  of 
$\fh$-weights occurring in $V$ is finite with 
$V = \bigoplus_{\beta\in \fh^*} V_\beta(\fh).$
\end{lemma}

\begin{proof} The derived representation 
$\L(\pi) \: \fh \to \gl(V)$ is a homomorphism of Banach--Lie algebras. 
For each character $\chi \in \hat T$ for which $V_\chi(T)$ is non-zero, 
we therefore have $\|\L(\chi)\| \leq \|\L(\pi)\|$, and since 
the character group 
$\hat T$ is a discrete subgroup of the dual $\L(T)^* \cong \fh_\R^*$, the set 
${\cal P}(V,\fh)$ is finite. The remaining assertion follows 
from Theorem~\ref{thm:peter-weyl}(3). 
\end{proof} 

\begin{definition} (The derived representation) \label{def:2.2}
Let $E$ be a locally convex space and 
$V := {\cal O}(G,E)$ denote the space of $E$-valued 
holomorphic functions on $G$ on which $G$ acts by 
$$ (\pi(g)f)(h) := f(g^{-1}h). $$
We endow $V$ with the compact open topology, turning it into a locally convex 
space. Then, for each $f \in V$,  the map 
$$ \pi^f \: G \to V, \quad g \mapsto \pi(g)f $$
is holomorphic because the map 
$$ \hat\pi^f \: 
G \times G \to E, \quad (g,h) \mapsto (\pi(g)f)(h) = f(g^{-1}h) $$
is holomorphic (\cite[Prop.~III.13]{Ne01}; \cite[Thm.~2.2.5]{Mue06}). 
We thus obtain an equivariant embedding 
$V \into {\cal O}(G,V), f \mapsto \pi^f$, 
and from that we derive that 
$$ \L(\pi)(x)f = - X_r f, $$
where $X_r \in {\cal V}(G)$ is the right invariant vector field 
with $X_r(\1) = x$, defines a derived representation of $\g$ on $V$. 
It is called the {\it derived representation $\L(\pi)$}. 
\end{definition}

\begin{definition} 
Let $P \cong L \ltimes U \subeq G$ be a 
parabolic subgroup (cf.\ Proposition~\ref{prop:1.9}) and 
$(\rho,W)$ be a holomorphic representation of the weakly root graded 
Banach--Lie group $L$ on $E$. We extend this 
representation in the canonical fashion to a representation 
$\rho$ of $P$ with $U \subeq \ker \rho$. 
We obtain a representation 
on the space 
$$ {\cal O}_\rho(G,E) 
:= \{ f \in {\cal O}(G,E) \: (\forall g \in G)(\forall p \in P)\ 
f(gp) = \rho(p)^{-1}f(g)\} $$
by $(\pi(g)f)(x) = f(g^{-1}x).$
\end{definition} 

\begin{lemma} \label{lem:2.5b} 
Let $(\rho,E)$ be any holomorphic representation of $P$. 
If we endow the space ${\cal O}_\rho(G,E)$ with the 
compact open topology, then the following assertions hold: 
\begin{description}
\item[\rm(a)] The action of $G_\Delta$ on the complex locally 
convex space ${\cal O}_\rho(G,E)$ is holomorphic. 
\item[\rm(b)] The space  ${\cal O}_\rho(G,E)^{[\fh]}$ of $\fh$-finite 
vectors is dense. 
\end{description}
\end{lemma}

\begin{proof} (a) is a direct consequence of Theorem~\ref{thm:action} in 
Appendix~B.

(b) First we use \cite[Thm.~III.11(c)]{Ne01} to see that 
the locally convex space ${\cal O}(G,E)$ is complete because
$G$ is Banach and $E$ is complete. In particular, the closed 
subspace ${\cal O}_\rho(G,E)$ is complete. 

In view of (a), the complex group $T_\C = H$ acts continuously 
on the complete locally convex space ${\cal O}_\rho(G,E)$, 
so that the assertion follows from the Peter--Weyl 
Theorem~\ref{thm:peter-weyl} and the fact that $T$-eigenvectors 
for a character $\chi$ are $\fh$-eigenvectors for $\L(\chi) \in \fh^*$. 
\end{proof}

\begin{defn} \label{def:taylor} 
We assign to each $D \in U(\g)$ the 
right invariant differential operator $D_r$ on $G$, such that 
$D \mapsto D_r$ is an anti-homomorphism of associative algebras mapping 
$x \in \g$ to $X_r$, the right invariant vector field with 
$X_r(\1) = x$. It acts on $\cO(G,E)$ by 
$$ (X_r f)(g) := \derat0 f(\exp_G(tx)g). $$

Since $G$ is connected, each holomorphic 
function $f \: G \to E$ is determined by its jet 
in the identity. We may represent this jet 
by a linear map 
$$ \tilde f \: U(\g) \to E, \quad D \mapsto (D_r f)(\1). $$
For $x \in \fp$ we then have 
\begin{align*}
\tilde f(D x) 
&= ((Dx)_r f)(\1) 
= (X_r (D_r f))(\1) 
= T_\1(D_r f)(x) \\
&= - \L(\rho)(x)(D_rf)(\1) = -\L(\rho)(x)\tilde f(D), 
\end{align*}
showing that we obtain an injection 
\begin{equation}
  \label{eq:taylorphi}
\Phi \: {\cal O}_\rho(G,E) \into \Hom_\fp(U(\g), E)_c, \quad 
\Phi(f)(D) := (D_r f)(\1)
\end{equation}
into the coinduced Lie algebra representation. The continuity 
of $\Phi(f)$ on $U(\g)$ follows from the Taylor expansion of the 
holomorphic function $f$ in~$\1$.
The $\g$-equivariance of this injection follows from 
\begin{align*}
\Phi(xf)(D) 
= -\big((D_r X_r)f\big)(\1) 
= -\big((xD)_r f\big)(\1) 
= - \Phi(f)(xD) = (x.\Phi(f))(D).
\end{align*}
\end{defn}

We now come to the main result of this section. Under more restrictive 
assumptions on $(\rho,E)$, the following theorem will be sharpened 
in Theorem~\ref{thm:extens} below. 

\begin{theorem} \label{thm:fin-weight} If $G$ is a connected weakly 
$\Delta$-graded Lie group, $P$ a parabolic subgroup and 
$(\rho,E)$ a holomorphic representation of $P$ which is 
trivial on $U$, then 
${\cal O}_\rho(G,E)$ is a finite direct sum of 
$\fh$-weight spaces. 

If, in addition, 
$P$ is connected, then 
${\cal O}_\rho(G,E)$ carries the structure of a holomorphic 
Banach $G$-module. If, moreover, $G$ is simply connected, 
$$ \Phi \: {\cal O}_\rho(G,E)\to \Hom_\fp(U(\g),E)_c^{[\g_\Delta]}, \quad 
\Phi(f)(D) := (D_r f)(\1) $$ 
is an isomorphism of $\g$-modules. 
\end{theorem}

\begin{proof} 
In view of Lemma~\ref{lem:2.1}, $E$ decomposes into finitely many $\h$-weight 
spaces, so that Theorem~\ref{thm:weightfin-alg} implies that the set 
$\cP(\cO_\rho(G,E),\fh)$ of $\fh$-weights in the 
submodule $\cO_\rho(G,E)^{[\g_\Delta]}$ 
of $\g_\Delta$-finite elements is finite. Note that, in view of 
Lemma~\ref{lem:2.5b} and the connectedness of $G_\Delta$, 
this space coincides with the space of 
$G_\Delta$-finite elements. 

Applying the Peter--Weyl Theorem to a maximal compact subgroup 
of $G_\Delta$ (Lemma~\ref{lem:2.5b}), 
we see that $\cO_\rho(G,E)^{[\g_\Delta]}$ 
is dense in $\cO_\rho(G,E)$ 
with respect to the compact open topology. 
Then we apply the Peter--Weyl Theorem~\ref{thm:peter-weyl} to the 
$T$-action on $\cO_\rho(G,E)$ to obtain that each $T$-, resp., 
$\h$-weight of $\cO_\rho(G,E)$ occurs in the dense subspace  
$\cO_\rho(G,E)^{[\g_\Delta]}$. 
We conclude that 
$\cO_\rho(G,E)$ is a finite direct sum of $\fh$-weight spaces. 

This in turn implies that $\cO_\rho(G,E)$ is a locally finite 
$\g_\Delta$-module \break (Lemma~\ref{lem:fin-iso}) and hence that 
$\im(\Phi)$ 
is contained in the locally $\fh$-finite subspace 
$V_c := \Hom_\fp(U(\g), E)_c^{[\fh]}$ 
of the coinduced 
representation $\Hom_\fp(U(\g), E)_c$.  
Since $\cO_\rho(G,E)$ is locally finite under $\g_\Delta$, 
the same holds for its image under $\Phi$, and thus 
$\im(\Phi) \subeq V_f := V_c^{[\g_\Delta]}$. 

Now we assume that $P$ is connected. 
We recall from Theorem~\ref{thm:coind-ban} that $V_f$ is a Banach 
$\g$-module. Let $q_G \: \tilde G \to G$ be the universal covering morphism,  
$\tilde P := q_G^{-1}(P)_0$ be the identity component of the 
inverse image of $P$ and $\tilde\rho := \rho \circ q_G\res_{\tilde P}$ 
be the representation of $\tilde P$ on $E$ obtained from $\rho$. 
Then the action of $\g$ on 
$V_f$ integrates to a unique holomorphic $\tilde G$-representation 
given by a morphism $\tilde G \to \GL(V_f)$ of Banach--Lie groups 
(\cite[Ch.~III, \S 6.1, Thm~1]{Bou89}. Since 
$\tilde G$ is connected, the corresponding map 
$\tilde\Phi \: \cO_{\tilde\rho}(\tilde G, E) \to V_f$
is $\tilde G$-equivariant (cf.\ \cite{GN09}). 
{}From Corollary~\ref{cor:5.5} in the appendix, we further derive that 
$V_f \subeq\im(\tilde\Phi)$, which proves surjectivity. 
Since the pullback map 
$$ q_G^* \: \cO_{\rho}(G,E) \to \cO_{\tilde\rho}(\tilde G,E), \quad 
f \mapsto f \circ q_G $$
is an injection whose range coincides with the set of functions  
constant on the cosets of $\ker q_G$, 
$\cO_{\rho}(G,E)$ 
can be identified with the closed subspace of 
$\ker q_G$-fixed points in the Banach $\tilde G$-module 
$V_f$, on which the representation of 
$\tilde G$ factors through a holomorphic 
representation of $G \cong \tilde G/\ker q_G$. This completes the proof. 
\end{proof}

If $G$ is a finite dimensional complex reductive group, then 
all homogeneous spaces $G/P$ are compact, so that it follows immediately 
from Liouville's Theorem that all holomorphic functions thereon are 
constant. If $G$ is infinite dimensional, then so is $G/P$.
In particular, it is never compact (locally compact Banach spaces 
are finite dimensional). However, it behaves in many respects 
like a compact flag manifold. Combining the preceding theorem with 
Remark~\ref{rem:triv-mod}, we shall see below that all holomorphic 
functions on any $G/P$ are constant. This result complements naturally 
the results of Dineen and Mellon on symmetric Banach manifolds of 
compact type \cite{DM98}. Such manifolds arise naturally 
from Banach--Lie groups graded by the root system $A_1$, and their 
Jordan theoretic proof essentially reduces matters to the case 
$G = \GL_2(C(X))$, where $X$ is a compact Hausdorff space. 
To put the following corollary into proper perspective, 
one should note that it also applies to all finite dimensional 
root graded groups which are not reductive. A typical example 
is $\SL_n(A)$ for any finite dimensional complex algebra~$A$. 
For $A = \C[\eps]$, $\eps^2 = 0$, we obtain in particular 
the tangent bundle $T(\bP_1(\C))$ as $\SL_2(A)/P$ where $P \subeq \SL_2(A)$ 
is the parabolic subgroup of upper triangular matrices 
(cf.\ \cite{NS09}).

\begin{cor} If $G$ is a connected root graded Banach--Lie group 
and $P \subeq G$ a parabolic subgroup, 
then all holomorphic functions on $G/P$ are constant. 
\end{cor} 

\begin{prf} If $q_G \: \tilde G \to G$ is the universal covering map 
and $\tilde P := q_G^{-1}(P)$, then $\tilde P$ is a parabolic subgroup 
of $\tilde G$ with 
$\tilde G/\tilde P \cong G/P$. We may therefore assume that 
$G$ is simply connected. 
If $P$ is not connected, then the canonical map $q \: G/P_0 \to G/P$ is 
a covering, and since each holomorphic function on $G/P$ pulls back 
to a holomorphic function on $G/P_0$, we may also assume that 
$P$ is connected. For the trivial representation $\rho \: P \to \C^\times \cong \GL(\C)$ we then
have 
$$ \cO(G/P) \cong \cO_\rho(G,\C) \cong \Hom_\fp(U(\g),\C)_c^{[\g_\Delta]} $$ 
(Theorem~\ref{thm:fin-weight}),  
and Remark~\ref{rem:triv-mod} implies that this space is 
one dimensional with trivial $G$-action. 
Therefore each holomorphic function on $G/P$ is constant. 
\end{prf}

\begin{corollary} \label{cor:findim} If $G$ is finite dimensional connected, 
and $(\rho,E)$ is a finite dimensional 
holomorphic representation of $P$, then $\dim {\cal O}_\rho(G,E) < \infty$. 
\end{corollary}

\begin{prf} Since $\cO_\rho(G,E)$ embeds into 
$\Hom_\fp(U(\g),E)_c^{[\g_\Delta]}$, it suffices to show that the latter 
space is finite dimensional. In the proof of Theorem~\ref{thm:coind-ban},  
we obtained the Banach structure on this space by embedding it into 
\break $\Hom(U(\fn),E)_c^k$ for some $k \in \N$. 
If $\fg$, and therefore $\fn$, is finite dimensional, then 
$\Hom(U(\fn),E)_c^k$ is finite dimensional for each $k \in \N$, and 
\break $\Hom_\fp(U(\g),E)_c^{[\g_\Delta]}$ inherits this property. 
\end{prf}

\section{Non-connected parabolic subgroups} \label{sec:4}

In Theorem~\ref{thm:fin-weight}, we have seen how to identify the space 
$\cO_\rho(G,E)$ with the $\g_\Delta$-finite submodule of the 
corresponding coinduced Lie algebra module, provided $G$ is $1$-connected 
and $P$ is connected. If $G$ is connected but not simply connected, 
then we may always consider the universal covering map 
$q_G \: \tilde G \to G$ and 
consider the parabolic subgroup $\tilde P := q_G^{-1}(P)$ in $\tilde G$. 
For $\tilde \rho := \rho \circ q_G\res_{\tilde P}$, the pullback map 
then induces a bijection 
$$q_G^* \: \cO_\rho(G,E) \to \cO_{\tilde \rho}(\tilde G,E), 
\quad f \mapsto f \circ q_G. $$
In fact, since $\ker q_G$ is contained in $\ker \tilde\rho$, 
each $f \in \cO_{\tilde \rho}(\tilde G,E)$ 
is constant on the cosets of $\ker q_G$, hence factors through a 
function on $G$, which clearly is contained in $\cO_\rho(G,E)$. 

This argument shows that the assumption of $G$ being simply connected 
is quite harmless. However, the connectedness of $P$ is a more tricky 
issue, and if $G$ is simply connected, then 
$\pi_0(P) \cong \pi_1(G/P)$ is non-trivial if $G/P$ is not simply 
connected (cf.\ Remark~\ref{rem:5.1}). This happens for many natural 
examples (see Example~\ref{exam:5.5}). 

Therefore it is desirable to understand the passage from a 
parabolic subgroup $P$ to its identity component $P_0$ and, putting 
$\rho_0 := \rho\res_{P_0}$, to understand the difference between 
$\cO_{\rho_0}(G,E)$ and its subspace $\cO_{\rho}(G,E)$. 
Of particular interest is the question whether the non-triviality of 
$\cO_{\rho_0}(G,E)$ implies that $\cO_{\rho}(G,E)$ is also non-trivial, 
resp., for which extension $\rho$ of $\rho_0$, 
the space $\cO_\rho(G,E)$ is non-trivial. 

\begin{remark}
  \label{rem:5.1} 
Since $G$ is connected, the long exact homotopy sequence of the 
$P$-principal bundle $G \to G/P$ provides an exact sequence 
$$ \pi_1(P) \to \pi_1(G) \to \pi_1(G/P) \to \pi_0(P) \to \1, $$
and if $G$ is $1$-connected, this leads to 
$\pi_1(G/P) \cong \pi_0(P).$
\end{remark}

\begin{rem} \label{rem:5.2} 
%(a) The group $\pi_0(P)$ acts on $\Hom(P,\C^\times)$ by 
%$$ ([p].\eta)(x) := \eta(p^{-1}xp), $$ 
%and if $\eta_0$ extends to a character of $P$, then $\eta_0$ is invariant 
%under this action. 
 If a holomorphic representation $\rho \: P_1 \to \GL(E)$ 
extends to a larger subgroup $P \subeq N_G(P_1)$ containing $P_1$, 
then, for each $p \in P$, the representation 
$$ (p.\rho)(x) := \rho(p^{-1}xp) $$
is equivalent to $\rho$. Since the equivalence class of $p.\rho$ 
only depends on the coset $[p] \in P/P_1$, we obtain an action of 
$P/P_1$ on the set of equivalence classes of holomorphic 
representations of $P$ on $E$. 
\end{rem} 

\begin{rem} \label{rem:vecbun-iso} Let $g \in N_G(P)$ 
and consider the corresponding right multiplication 
$\oline\rho_g \: G/P \to G/P, xP \mapsto xgP$. 

For a holomorphic representation $\rho \: P \to \GL(E)$ we define 
the associated holomorphic homogeneous 
vector bundle 
$$\bE _\rho:= (G \times E)/P := G \times_\rho E$$ 
whose elements we write as $[g,v]$. 

(a) For two holomorphic representations 
$\rho_1, \rho_2 \: P \to \GL(E)$, any $G$-equivariant bundle isomorphism 
$\phi \: \bE_{\rho_1} \to \bE_{\rho_2}$ covering $\oline\rho_g$ 
is of the form $\phi([y,v]) = [yg, \psi(v)]$ for some 
$\psi \in \GL(E)$. Since $\phi$ is well-defined, we have for 
each $p \in P$: 
\begin{align*}
[yg, \psi(v)] = \phi([y,v]) 
&= \phi([yp, \rho_1(p)^{-1}v]) 
= [ypg, \psi\rho_1(p)^{-1}v] \\
&= [yg, \rho_2(g^{-1}pg)\psi\rho_1(p)^{-1}v], 
\end{align*}
which means that 
\begin{equation}
  \label{eq:intertwine}
\psi \circ \rho_1(p) = \rho_2(g^{-1}pg) \circ \psi \quad \mbox{ for } \quad 
p \in P, 
\end{equation}
i.e., that $\psi$ intertwines the representations $\rho_1$ and $g.\rho_2$. 
If, conversely, $\psi$ satisfies this relation, then 
$\phi([y,v]) := [yg, \psi(v)]$ is a well-defined $G$-equivariant 
bundle morphism covering $\oline\rho_g$. 

(b) For $g \in N_G(P)$ and $\psi \in \GL(E)$ we consider the operator 
$$ M_g \: \cO(G,E) \to \cO(G,E), \quad 
M_g(f) := \psi \circ f \circ \rho_g, \quad \rho_g(x) = xg. $$
If $M_g$ maps ${\cal O}_{\rho_1}(G,E)$ into $\cO_{\rho_2}(G,E)$, 
then 
$$  \rho_2(p)^{-1} \psi f(yg) = \psi f(ypg) 
= \psi f(ygg^{-1}pg) = \psi \rho_1(g^{-1}pg)^{-1} f(yg) $$
holds for each $f \in {\cal O}_{\rho_1}(G,E)$, $y \in G$ and $p \in P$. 
If this is the case and $\ev_\1\res_{\cO_{\rho_1}(G,E)}$ 
has dense range in $E$, then we obtain 
$\rho_2(p)^{-1} \psi  = \psi \rho_1(g^{-1}pg)^{-1},$
which is equivalent to  
$$  \rho_1(g^{-1}pg)  = \psi^{-1} \rho_2(p) \psi \quad \mbox{ for } \quad 
p \in P. $$
If, conversely, this condition is satisfied, 
then the preceding 
calculation shows that $M_g$ maps ${\cal O}_{\rho_1}(G,E)$ 
into $\cO_{\rho_2}(G,E)$. 
\end{rem}

\begin{prop} \label{prop:4.3} 
Suppose that $\ev_\1 \: {\cal O}_\rho(G,E) \to E$ 
has dense range and let $p \in N_G(P)$. 
Then the following are equivalent: 
\begin{description}
\item[\rm(a)] The two representations $\rho$ and $p.\rho$ define 
$G$-equivalent holomorphic vector bundles $\bE_\rho$ and 
$\bE_{p.\rho}$. 
\item[\rm(b)] The representations $\rho$ and $p.\rho$ are equivalent. 
\item[\rm(c)] The right multiplication map $\oline\rho_p \: G/P \to G/P, 
xP \mapsto xgP$ lifts to a $G$-equivariant bundle isomorphism 
$M_p \: \bE_\rho \to \bE_\rho$. 
\item[\rm(d)] There exists an $A_p \in \GL(E)$ for which the operator 
\begin{equation}
  \label{eq:inter}
M_p(f) := A_p \circ f \circ \rho_p, \quad \rho_p(x) = xp, 
\end{equation}
preserves ${\cal O}_\rho(G,E)$. 
\end{description}
\end{prop}

\begin{proof} (a) $\Longleftrightarrow$ (b) follows from 
Remark~\ref{rem:vecbun-iso}(a), applied with $g = \1$, 
$\rho_1 = \rho$ and $\rho_2 = p.\rho$.  

(b) $\Longleftrightarrow$ (c) follows from 
Remark~\ref{rem:vecbun-iso}(a), applied to 
$\rho_1 = \rho_2 = \rho$. 

(b) $\Longleftrightarrow$ (d) follows from 
Remark~\ref{rem:vecbun-iso}(d). 
\end{proof}

\begin{corollary} \label{cor:4.4} 
If $\chi \: P \to \C^\times$ is a holomorphic 
character,  ${\cal O}_\rho(G) \not=\{0\}$ and $g \in N_G(P)$, 
then the space ${\cal O}_\rho(G)$ is invariant 
under $\rho_g$ if and only if 
$g.\rho =  \rho$. 
\end{corollary}

\begin{thm}\label{thm:extens} 
Assume that $G$ is a connected weakly root graded Banach--Lie group, 
$P_0 \subeq G$ is a connected parabolic subgroup, and 
the holomorphic representation 
$\rho_0 \: P_0 \to \GL(E)$ is irreducible with 
$\End_{P_0}(E) = \C \1$ and that 
$\cO_{\rho_0}(G,E)\not=\{0\}$. We further assume that 
$P \subeq N_G(P_0) = N_G(\fp)$ is an open subgroup satisfying 
$p.\rho_0 \sim \rho_0$ for each $p \in P$. 
Then the following assertions hold: 
\begin{description}
\item[\rm(a)] There exists a unique extension 
$\rho \: P \to \GL(E)$ of $\rho_0$ satisfying 
$\cO_{\rho_0}(G,E) = \cO_{\rho}(G,E).$
For all other homomorphic extensions $\gamma \: P \to \GL(E)$, the space 
$\cO_\gamma(G,E)$ is trivial. 
\item[\rm(b)] The map 
$\Phi \: \cO_\rho(G,E) \to \Hom_\fp(U(\g),E)_c^{[\g_\Delta]}, 
\Phi(f)(D) = (D_r f)(\1)$
is a linear isomorphism.  
\end{description}
\end{thm}

\begin{prf} {\bf Case 1: $G$ is simply connected.} 
First we use Lemma~\ref{lem:2.1} to see that 
the representation $\rho_0$ has only finitely many $\fh$-weights. 
Since $G$ is simply connected, Theorem~\ref{thm:fin-weight} provides 
a $\g$-equivariant isomorphism 
$$ \Phi \: \cO_{\rho_0}(G,E) \to \Hom_\fp(U(\g),E)_c^{[\g_\Delta]}. $$
Therefore Proposition~\ref{prop:irred}(ii) implies that 
$\ev_\1 \: \cO_{\rho_0}(G,E)^N \to E$
is a linear isomorphism of the subspace of $N$-invariant elements onto $E$. 
In particular, $\ev_\1 \: \cO_{\rho_0}(G,E) \to E$ is surjective. 
For each $\phi \in \End_G(\cO_{\rho_0}(G,E)),$ we now derive from 
Proposition~\ref{prop:irred}(iii),(iv) that the $P_0$-morphism 
$$\ev_\1 \circ \phi \: \cO_{\rho_0}(G,E) \to E$$ 
annihilates $\ker(\ev_\1)$ 
and that there exists a $\psi = \alpha \id_E \in \End_\fp(E) = \C \1$ with 
$$ \ev_\1 \circ \phi = \psi \circ  \ev_\1 = \alpha \cdot \ev_\1. $$
This proves that, for each $g \in G$, we have  
$$ \phi(f)(g) = \ev_\1(g^{-1}.\phi(f)) 
= \ev_\1(\phi(g^{-1}.f)) 
= \alpha \ev_\1(g^{-1}.f) = \alpha f(g), $$
i.e., $\phi(f) = \alpha f$. We thus obtain 
\begin{equation}
  \label{eq:commutant}
\End_G(\cO_{\rho_0}(G,E)) = \C \id. 
\end{equation}

In the following we write $c_g(x) := gxg^{-1}$ for conjugation with~$g$. 
Now let $p \in P$, pick $A_p \in \GL(E)$ with 
$p.\rho_0 = 
c_{A_p} \circ \rho_0$, and define the operator 
$$ M_p \: \cO_{\rho_0}(G,E) \to \cO_{\rho_0}(G,E), \quad 
M_p f := A_p \circ f \circ \rho_p $$
(Proposition~\ref{prop:4.3}(d)). Since $M_p$ commutes with the 
$G$-action by left translations, there exists an $\alpha_p \in \C$ 
with $M_p = \alpha_p \1$. This means that each $f \in \cO_{\rho_0}(G,E)$ 
satisfies the equation 
$f \circ \rho_p = \alpha_p A_p^{-1} \circ f.$
Therefore 
$\rho(p) := \alpha_p^{-1} A_p \in \GL(E)$
satisfies 
\begin{equation}
  \label{eq:extrel}
f(xp) = \rho(p)^{-1} f(x) \quad \mbox{ for all } \quad x \in G, 
f \in \cO_{\rho_0}(G,E). 
\end{equation}
Since $\ev_\1$ is surjective, evaluating in $x  = \1$ shows that this 
relation determines $\rho(p)$ uniquely. In particular, 
$\rho(p) = \rho_0(p)$ for $p \in P_0$. 
We thus obtain a map 
$\rho \: P \to \GL(E)$, determined uniquely by \eqref{eq:extrel}. 
For $p_1, p_2 \in P$ we have for each $f \in \cO_{\rho_0}(G,E)$ 
$$ \rho(p_1p_2)^{-1} f(x) = f(xp_1 p_2)
= \rho(p_2)^{-1} \rho(p_1)^{-1} f(x), $$
showing that $\rho$ is multiplicative, hence a representation of 
$P$ on $E$. Clearly, our construction implies that 
$\cO_\rho(G,E) = \cO_{\rho_0}(G,E)$. 

If $\gamma \: P \to \GL(E)$ is another extension of 
$\rho_0 \: P_0\to \GL(E)$, then we also have 
$\rho_0 \circ c_p = c_{\gamma(p)} \circ \rho_0$
for each $p \in P$, so that the argument in the preceding 
proof implies that 
$\gamma = \chi \cdot \rho$ for some map 
$\chi \: P \to \C^\times$. Since $\chi(P)$ is central in $\GL(E)$, 
the fact that $\gamma$ is a homomorphism implies that 
$\chi$ is a homomorphism, and since $\gamma$ also 
extends $\rho_0$, $\chi$ vanishes on $P_0$. Assume that 
$\chi(p)\not=1$ for some $p \in P$. 

For any $f \in \cO_\gamma(G,E)$, we now have 
$f \in \cO_{\rho_0}(G,E) = \cO_\rho(G,E)$, and therefore 
$$ \chi(p)\rho(p) f(\1) = \gamma(p) f(\1) = f(p^{-1}) = \rho(p) f(\1). $$
Since $\chi(p)\not=1$, this implies that $f(\1) =0$, and since 
$\cO_\gamma(G,E)$ is $G$-invariant, we obtain $f = 0$. 

{\bf Case 2: $G$ is connected but not simply connected.} 
Let $q_G \: \tilde G \to G$ be the universal covering, 
$\hat P_1 := q_G^{-1}(P_0)$ and  $\hat P := q_G^{-1}(P)$. 
We also put 
$\hat \rho_1 := \rho_0 \circ q_G\res_{\hat P_1}$ and 
$\hat \rho_0 := \hat\rho_1 \res_{\hat P_0}$. 
Then we have a linear isomorphism 
$$ q_G^* \: \cO_{\rho_0}(G,E) \to \cO_{\hat\rho_1}(\tilde G,E), \quad 
f \mapsto f \circ q_G $$
because $\ker q_G \subeq \hat P_1$ is contained in the kernel of 
$\hat\rho_1$. From the simply connected case we derive that 
$\cO_{\hat\rho_1}(\tilde G,E) = \cO_{\hat\rho_0}(\tilde G,E)$. 

Since $\hat P$ satisfies $p.\hat \rho_0 \sim \hat\rho_0$ for each 
$p \in \hat P$, there  exists an extension 
$\hat\rho \: \hat P \to \GL(E)$ with 
$$ \cO_{\hat\rho}(\tilde G,E) =  \cO_{\hat\rho_0}(\tilde G,E). $$
{}From the uniqueness assertion in Case 1, we derive that 
$\hat\rho\res_{\hat P_1} = \hat\rho_1$, hence that 
$\ker q_G \subeq \ker\hat\rho$. Therefore 
$\hat\rho$ factors through an extension 
$\rho \: P \to \GL(E)$ with $\rho \circ q_G = \hat\rho$, for which 
$$ q_G^* \: \cO_{\rho}(G,E) \to 
\cO_{\hat\rho}(\tilde G,E)  
= \cO_{\hat\rho_0}(\tilde G,E) = q_G^*(\cO_{\rho_0}(G,E)) $$
is a linear isomorphism. This implies that 
$\cO_{\rho}(G,E) =  \cO_{\rho_0}(G,E)$
because $q_G^*$ is injective. For any other extension 
$\gamma$ of $\rho_0$ we obtain an extension 
$\hat\gamma := \gamma \circ q_G$ of $\hat\rho_0$ to $\hat P$ which differs 
from $\rho$. Therefore 
$q_G^*(\cO_{\gamma}(G,E)) \subeq 
\cO_{\hat\gamma}(\tilde G,E) = \{0\}$ implies that 
$\cO_\gamma(G,E)$ vanishes. 

This completes the proof of (a). To verify (b), we collect the information 
obtained so far to obtain isomorphisms 
$$ \cO_\rho(G,E) \ssmapright{q_G^*} \cO_{\hat\rho}(\tilde G,E) 
= \cO_{\hat\rho_0}(\tilde G,E) 
\ssmapright{\Phi_{\tilde G}} \Hom_\fp(U(\g),E)_c^{[\g_\Delta]}. $$
Since $\Phi_{\tilde G} \circ q_G^* = \Phi$, this proves (b). 
\end{prf}

\begin{cor}\label{thm:extens-line} 
Assume that $G$ is connected and that 
the holomorphic character 
$\chi_0 \: P_0 \to \C^\times$ satisfies $\cO_{\chi_0}(G)\not=\{0\}$. 
We further assume that 
$P \subeq N_G(P_0)$ is an open subgroup fixing $\chi_0$. 
Then there exists a unique extension 
$\chi \: P \to \C^\times$ of $\rho_0$ to a holomorphic character of 
$P$ such that 
$$ \cO_{\chi_0}(G) = \cO_{\chi}(G). $$
For all other extensions $\gamma$ of $\chi_0$, we have 
$\cO_\gamma(G) = \{0\}$. 
\end{cor}

\begin{rem} If $G$ is finite dimensional and complex reductive, 
i.e., $Z(G)_0 \cong (\C^\times)^r$ is a complex torus, and 
$E$ is also finite-dimensional, then $\cO_{\rho_0}(G,E)$ 
is finite dimensional by Corollary~\ref{cor:findim}. 
This implies that the holomorphic representation of $G$ on 
this space is semisimple, so that the relation 
\eqref{eq:commutant} from the proof of Theorem~\ref{thm:extens} 
implies that $\cO_{\rho_0}(G,E)$ is an irreducible $G$-representation. 
\end{rem}

\begin{example}
  \label{exam:5.5}
For $A = C(\bS^1, \C)$ and $G = \SL_2(A)_0$ (the identity component), 
the standard 
parabolic subgroup $P$ is given by 
$$P = \Big\{ 
\begin{pmatrix} a & b \cr 0 & a^{-1}\end{pmatrix} \: 
a \in A^\times, b \in A \Big\} 
\cong A \rtimes A^\times,$$ 
because for each $a \in A^\times$, the relation 
$$ \pmat{a & 0 \\ 0 & a^{-1}}  
= \pmat{1 & 0 \\ a^{-1}-1 & 1}  
\pmat{1 & 1 \\ 0 & 1}  
\pmat{1 & 0 \\ a-1 & 1}  
\pmat{1 & -a^{-1} \\ 0& 1}  $$
implies that all diagonal matrices in $\SL_2(A)$ are lying in the 
identity component. We conclude that 
$\pi_0(P) \cong \pi_0(A^\times) \cong \Z$. 
As $A^\times$ is abelian and $A^\times_0$ is divisible, we have 
$$ A^\times \cong A^\times_0 \times \Z 
\quad \mbox{ and } \quad P \cong P_0 \rtimes \Z. $$
In view of 
$(P_0,P_0) = A \rtimes \{\1\}$, each character of $P_0$
factors through $A^\times_0$. We also see that each character 
has many extensions to a character on $P$. 
In view of Theorem~\ref{thm:extens}, 
this implies the existence of extensions 
$\chi$ with 
  $$ \{0\} = {\cal O}_\chi(G) \not= {\cal O}_{\chi_0}(G).$$

Typical examples arise as follows. If $\chi(a) := a(x)$ arises from 
evaluation in some $x \in \bS^1$, then 
${\cal O}_\chi(G)$ is non-zero, and the corresponding space of 
holomorphic sections is isomorphic to $\C^2$, the representation 
on this space being given by 
$$ G \to \SL_2(\C), \quad g \mapsto g(x)^{-\top}. $$

The winding number 
$$ w \: A^\times \to \Z $$
is a character of $A^\times$, vanishing on $A^\times_0$, for which 
all other 
extensions of $\chi_0$ to $A^\times$ are of the form 
$$ \chi_z(a) = \chi(a)z^{w(a)} $$
for some $z \in \C^\times$. 
In view of Theorem~\ref{thm:extens}, we obtain only for 
$z = 1$ a non-trivial space of holomorphic sections. 
\end{example}

\section{Realizing irreducible holomorphic representations} \label{sec:5}

In the preceding sections we developed techniques to study the 
holomorphic representations of $G$ in the Banach spaces 
$\cO_\rho(G,E)$, defined by a holomorphic representation 
$(\rho,E)$ of $P$. 
In this section we briefly discuss the converse, namely to which 
extent irreducible holomorphic Banach representations $(V,\pi)$ of the connected
group $G$ have realizations in spaces $\cO_\rho(G,P)$, where $U \subeq \ker \rho$. 

Let $(V,\pi)$ be an irreducible holomorphic Banach representation of 
$G$, i.e., all closed $G$-invariant subspaces of $V$ are trivial. 
Since the representation is holomorphic, each closed $\g$-invariant 
subspace is $G$-invariant, so that the derived $\g$-representation on~$V$ 
is also irreducible. We consider a connected parabolic subgroup 
$P$ of~$G$. % and assume that $G$ is simply connected.

We pick $x_\Sigma$ as in Remark~\ref{rem:para-pos} and recall from 
Lemma~\ref{lem:2.1} that $V$ is the direct sum of finitely many $\fh$-weight 
spaces. Let $\lambda \in \R$ denote the minimal eigenvalue of 
$x_\Sigma$ on $V$. Since $\Sigma^-(x_\Sigma) \leq -1$, we clearly have 
\begin{equation}
  \label{eq:5.1}
V_\lambda(x_\Sigma) \subeq V^\fn, 
\end{equation}
and likewise 
\begin{equation}
  \label{eq:5.2}
  V_{>\lambda}(x_\Sigma) := \sum_{\mu > \lambda} V_\mu(x_\Sigma) \supeq \fu.V.
\end{equation}
Since $x_\Sigma$ is central in $\fl$, the closed subspace 
$V_\lambda(x_\Sigma)$ is $\fl$-invariant. For each closed $\fl$-invariant subspace 
$W \subeq V_\lambda(x_\Sigma) \subeq V^\fn$, the PBW Theorem implies that 
$$ U(\g)W = U(\fu)U(\fl)U(\fn) W 
= U(\fu)W \subeq W + \fu.V \subeq W \oplus V_{> \lambda}(x_\Sigma). $$
Therefore the irreducibility of the $\g$-module $V$ implies that 
$V_\lambda(x_\Sigma)$ is an irreducible $\fl$-module. Further, the density of 
$U(\g)V_\lambda(x_\Sigma)$ in $V$ implies the equality 
\begin{equation}
  \label{eq:5.3}
\oline{\fu.V} = V_{> \lambda}(x_\Sigma). 
\end{equation}

Since $\oline{\fu.V}$ is $\fp$-invariant, the quotient space 
$E := V/\oline{\fu.V}$ 
inherits a natural $\fp$-module structure, where $\fu$ acts trivially. 
We write 
$\beta \: V \to E$
for the corresponding quotient map. Let $\rho$ denote the corresponding 
representation of $P$ on this Banach space. Then Theorem~\ref{thm:a.5}(2) 
yields an embedding 
$$ \beta_G \: V \into \cO_{\rho}(G,E), \quad 
\beta_G(v)(g) := \beta(g^{-1}.v). $$
Clearly, $\beta_G$ maps $V_\lambda(x_\Sigma)$ 
into the space of $N$-invariants 
in $\cO_{\rho}(G,E)$ and satisfies $\ev_\1 \circ \beta_G = \beta$, 
so that $\ev_1 \: \cO_{\rho}(G,E) \to E$ is surjective. 
Therefore 
$$ \beta_G(V_\lambda(x_\Sigma)) = \cO_{\rho}(G,E)^N \cong E. $$ 

To see that $\beta_G$ is continuous with respect to the natural 
Banach structure on $\cO_\rho(G,E)$, we observe that for 
$x_1, \ldots, x_k \in \g$, we have 
$$ \Phi(\beta_G(v))(x_1 \cdots x_k) 
= ((X_k)_r \cdots (X_1)_r \beta_G(v))(\1) 
= (-1)^k \beta(x_k \cdots x_1.v), $$
which defined a continuous $k$-linear map $\g^k \to E$. 
{}From the continuity of $\beta_G$ we finally derive that 
its image is contained in the closed $G$-submodule of 
$\cO_\rho(G,E)$ generated by the subspace 
$\cO_\rho(G,E)^N \cong E$ of $N$-invariants. 
(Proposition~\ref{prop:irred}(ii)). 

%and the image of $V$ is dense in the minimal $\g$-submodule of 
%$\cO_\rho(G,E) \cong \Hom_\fp(U(\g),E)_c^{[\g_\Delta]}$ 
%(Proposition~\ref{prop:irred} and Theorem~\ref{thm:fin-weight}; 
%here we use that $G$ is simply connected). 

Putting all this together, we have: 
\begin{thm}\label{thm:realize} 
Let $G$ be a connected weakly root graded Banach--Lie group 
and $(\pi,V)$ be an irreducible holomorphic representation of~$G$. 
If $P$ is a connected parabolic subgroup of $G$, then 
we obtain an irreducible holomorphic representation 
$(\rho,E)$ on $E := V/\oline{\fu.V}$ with $U \subeq \ker \rho$. 
If $\beta \: V \to E$ denotes the quotient map, then 
$$ \beta_G \: V \to \cO_\rho(G,E), \quad 
\beta_G(v)(g) := \beta(g^{-1}.v) $$
defines a continuous morphism of holomorphic Banach $G$-modules 
whose image is a dense subspace in the minimal closed 
$G$-submodule of $\cO_\rho(G,E)$. 
\end{thm}

\appendix
\section{Generalities on holomorphic induction} \label{sec:app}

The following theorem is the natural version of Frobenius reciprocity 
for holomorphically induced representations. The assumptions on the 
$G$-module $W$ and the $P$-module $E$ are rather weak, so that they 
do not only apply to Banach representations, but also to representations 
in spaces of holomorphic functions. 

\begin{theorem}{\rm(Frobenius Reciprocity)}  \label{thm:frob}
Let $W$ be a $G$-representation with holomorphic orbit maps and 
$(\rho,E)$ a locally convex representation of $P$ 
with holomorphic orbit maps. Then the map 
$$ \ev_\1 \circ \: \Hom_G(W, {\cal O}_\rho(G,E)) \to \Hom_P(W,E) $$
is bijective if $\Hom_P$ denotes the set of continuous $P$-morphisms 
and $\Hom_G$ the set of those 
$G$-morphisms that are continuous with respect to the 
topology of pointwise convergence on ${\cal O}_\rho(G,E)$. 
Its inverse is given by 
$$ \Hom_P(W,E) \to \Hom_G(W, \cO_\rho(G,E)), \quad 
\beta \mapsto \beta_G, \quad 
\beta_G(w)(g) = \beta(g^{-1}.w). $$
\end{theorem}

\begin{proof} Note that 
$$\ev_\1 \:  {\cal O}_\rho(G,E) \to E, \quad 
f \mapsto f(\1) $$
is $P$-equivariant: 
$$ (p.f)(\1) = f(p^{-1}) = \rho(p).f(\1) $$ 
and continuous with respect to the topology of pointwise convergence, 
so that $\ev_\1 \circ$ maps $\Hom_G(W,\cO_\rho(G,E))$ into 
$\Hom_P(W,E)$. 

Let $\phi \in \Hom_G(W,\cO_\rho(G,E))$ and assume that 
$\ev_\1 \circ \phi = 0$, i.e., 
$$\phi(w)(\1) = 0 \quad \mbox{ for each } w \in W.$$ 
Since $\phi$ is $G$-equivariant, 
$\phi(w)$ vanishes on all of $G$, so that 
$\phi = 0$. Therefore $\ev_\1 \circ$ is injective. 

To see that it is also surjective, let $\beta \in \Hom_P(W,E)$ and define  
$$ \beta_G \: W \to {\cal O}(G,E), \quad 
\beta_G(w)(g) := \beta(g^{-1}.w) $$
(recall that $G$-orbit maps in $W$ are holomorphic).  
Then 
$$ \beta_G(w)(gp) = \beta(p^{-1}g^{-1}w) = \rho(p)^{-1}\beta(g^{-1}w) $$
implies that $\beta_G(W) \subeq {\cal O}_\rho(G,E)$, and it is clear that 
$\beta_G$ is $G$-equivariant with $\ev_\1 \circ \beta_G = \beta$. 
Finally we note that $\beta_G$ is continuous with respect to the topology 
of pointwise convergence on ${\cal O}(G,E)$ since $\beta$ is continuous.  
\end{proof} 

\begin{remark} (a) Note that the orbit maps for the $G$-action 
on the space ${\cal O}_\rho(G,E)$, endowed with the compact open 
topology, are holomorphic because for each $f \in {\cal O}_\rho(G,E)$, 
the map 
$$ G \times G \to E, \quad (x,y) \mapsto f(x^{-1}y) $$
is holomorphic (\cite[Prop.~III.13]{Ne01}). 

(b) In view of (a), each $G$-submodule $W$ of ${\cal O}_\rho(G,E)$ 
satisfies the assumptions of Theorem~\ref{thm:frob}, 
so that all $G$-morphisms $W \to {\cal O}_\rho(G,E)$, continuous 
w.r.t.\ the topology of pointwise convergence, 
correspond to continuous $P$-morphisms 
$W \to E$. 
\end{remark}

\begin{remark}   \label{rem:4.2} 
Let $(\pi,V)$ he a representation of $G$ with holomorphic orbit maps 
and 
$${\cal O}_\pi(G,V) := \{ f \in {\cal O}(G,V) \: 
(\forall g,x \in G)\ f(xg) = \pi(g)^{-1} f(x)\}, $$
then the evaluation map 
$$ \ev_\1 \: {\cal O}_\pi(G,V) \to V, \quad f \mapsto f(\1) $$
is a $G$-equivariant isomorphism whose inverse is given by 
$v \mapsto f_v$ with $f_v(g) = \pi(g^{-1})v$. 
Therefore we may identify ${\cal O}_\pi(G,V)$ with $V$. 
If we endow ${\cal O}_\pi(G,V)$ wit the topology of pointwise or 
compact convergence, we even obtain a topological isomorphism 
$V \cong {\cal O}_\pi(G,V)$. 
\end{remark}

\begin{cor}
${\cal O}_\rho(G,E)$ is non-zero if and only 
there exists some holomorphic $G$-module $W$ and a non-zero 
continuous $P$-homomorphism $W \to E$. 
\end{cor}

\begin{ex} Let $A$ be a unital Banach algebra and consider the 
root graded Banach--Lie group $G := \GL_n(A)$. 
Then the space 
$$ \fp = \{ (x_{ij}) \in \gl_n(A) \: i > j \Rarrow x_{ij}=0 \} $$
of upper triangular matrices is a parabolic subalgebra and 
$$ P := \fp \cap \GL_n(A) $$
is a corresponding parabolic subgroup. 

The preceding corollary provides a rich supply of holomorphic 
representations of $(\rho, E)$ of $P$ for which 
$\cO_\rho(G,E)$ is non-trivial. 
In the holomorphic $G$-module $W = A^n$, the subspaces 
$$ F_i := \sum_{j \leq i} A e_j, \quad e_j = (\delta_{ij})_{i=1,\ldots, n},  $$
are $P$-invariant, so that each space 
$E_i := A^n/F_i$ carries a natural holomorphic $P$-representation 
$\rho_i$ for which $\cO_{\rho_i}(G,E_i)\not=\{0\}$. 
Identifying $E_i$ in the natural fashion with $A^{n-i}$, the representation 
$\rho_i$ is given by 
$\rho(g_{kl}) = (g_{kl})_{i+1 \leq k,l \leq n},$
resp., 
$$ \rho\pmat{a & b \\ 0 & d} =  d, $$
if we write $X \in M_n(A)$ as a block matrix with entries  in 
$M_{i,i}(A)$, $M_{i,n-i}(A)$, $M_{n-i,i}(A)$ and $M_{n-i,n-i}(A)$, respectively. \end{ex}

\begin{theorem} \label{thm:a.5} Let $P \subeq G$ be a connected complex Lie 
subgroup, $(\pi,V)$ a 
$G$-representation with holomorphic orbit maps and 
$(\rho,E)$ a $P$-representation with holomorphic orbit maps. 
Then the following are equivalent: 
\begin{description}
\item[\rm(1)] There exists a $G$-cyclic continuous $\beta \in \Hom_P(V,E)$, 
i.e., $\beta(G.v) = \{0\}$ implies $v = 0$. 
\item[\rm(2)] There is a $G$-equivariant injection $\beta_G \: V \into 
{\cal O}_\rho(G,E)$ which is continuous with respect to the pointwise 
topology on ${\cal O}_\rho(G,E)$. 
\item[\rm(3)] $V$ embeds into $\Hom_\fp(U(\g),E)$ such that 
$\ev_\1 \: V \to E$ is continuous. 
\end{description}
\end{theorem}

\begin{proof} The equivalence of (1) and (2) follows from Theorem~\ref{thm:frob}. 

(2) $\Rarrow$ (3): For this implication, 
we only have to recall the inclusion 
$$ \Phi \: {\cal O}_\rho(G,W) \into \Hom_\fp(U(\g),E), \quad 
\Phi(f)(D) := (D_r f)(\1) $$
from Definition~\ref{def:taylor}. 

(3) $\Rarrow$ (1): The continuous linear map 
$\beta := \ev_\1\res_V \: V \to W$ is $\fp$-equivariant because  
$$ \ev_\1(p.\phi) = -\phi(p) = \L(\rho)(p)\phi(\1)
= \L(\rho)(p)\ev_\1(\phi), $$
hence $P$-equivariant because $P$ is connected (cf.\ \cite{GN09}). 
Further, $\beta$  is $\g$-cyclic because 
$0 = \beta(U(\g).f) = f(U(\g))$ 
implies $f = 0$. Since the $G$-orbit maps 
in $V$ are holomorphic, $\beta$ is also $G$-cyclic. 
\end{proof}

\begin{cor} \label{cor:5.5} Suppose that $P$ is connected and that 
$(\rho,E)$ is a $P$-rep\-re\-sen\-tation with holomorphic orbit maps. 
Then the image of the Taylor series map 
$$ \Phi  \: {\cal O}_\rho(G,E) \to \Hom_\fp(U(\g),E), \quad 
\Phi(f)(D) := (D_r f)(\1) $$
is the largest $\g$-submodule 
$V$ of $\Hom_\fp(U(\g),E)$ on which 
the $\g$-module structure integrates to a $G$-representation 
such that, for each $v \in V$, the map 
$$ G \to E, \quad g \mapsto (g.v)(\1) $$
is holomorphic. 
\end{cor}

\begin{proof} We recall from Definition~\ref{def:taylor} that 
$\Phi$ is injective and $\g$-equivariant. 
Let $V \subeq \Hom_\fp(U(\g),E)$ be a $\g$-submodule 
on which the representation integrates to a $G$-representation 
with the required properties. For each $v \in V$, 
we then obtain a holomorphic function 
$f_v \in {\cal O}(G,E)$ by $f_v(g) := (g^{-1}.v)(\1)$. 
Then $X_r f_v = -f_{x.v}$ holds for each $x \in \g$, and therefore 
$$ D_r f_v = f_{D^\sigma.v}, $$
where $\sigma \: U(\g) \to U(\g), D \mapsto D^\sigma$, 
is the unique antiautomorphism with $x^\sigma = - x$ for $x \in \g$. 
We now obtain 
$$ T(f_v)(D) = (D_r f_v)(\1) 
=  f_{D^\sigma.v}(\1) 
= (D^\sigma.v)(\1) = v(D), $$
so that $T(f_v) = v$. 
If $v \in V$, then $G.v \subeq \Hom_\fp(U(\g),E)$, so that we have 
for each $x \in \fp$  
$$ (X_l f_v)(g) 
= (- x.(g^{-1}.v))(\1) 
= (g^{-1}.v)(x) 
= - x.\big((g^{-1}.v)(\1)\big)
= - x.f_v(g), $$
and since $P$ is connected, we obtain $f_v \in {\cal O}_\rho(G,E)$ 
(cf.\ \cite{GN09}). 
\end{proof}

\section{A general continuity lemma} \label{sec:app2}

\begin{lemma} \label{lem:action} 
Let $M$ be a Hausdorff space, $V$ a locally convex
space, and $S$ a locally compact topological
semigroup which acts continuously on $M$ from the right. Then the
induced action
\begin{equation*}
\Phi\: S\times C(M,V)_{c}\to C(M,V)_{c},\quad 
(s,f)\nobreak\mapsto\nobreak(x\nobreak\mapsto\nobreak f(x.s))
\end{equation*}
is continuous with respect to the compact-open topology on the
function space.
\end{lemma}

\begin{proof} (\cite[Lemma~2.2.6]{Mue06})
Consider $s\in S$ and $f\in C(M,V)_{c}$. The basic open neighborhoods
of the image $\varphi:=\Phi(s,f)$ are of the form 
\begin{equation*}
U_{p,K,\varepsilon}(\varphi)=\{g\in C(M,V) \mid p(g(x)-\varphi(x))<\varepsilon
\text{ for all }x\in K\},
\end{equation*}
where $p$ is a continuous seminorm on $V$, $K\subseteq M$ is a compact
subset, and $\varepsilon>0$. We need to find an open neighborhood of
$(s,f)$ in $S\times C(M,V)$ that is mapped completely into
$U_{p,K,\varepsilon}(\varphi)$ by $\Phi$. To do this, we look at the
following inequality for $t\in S$ and $g\in C(M,V)_{c}$:
\begin{equation*}
p((t.g)(x)-\varphi(x))\leq p(g(x.t)-f(x.t))+p(f(x.t)-f(x.s))
\end{equation*}
By the continuity of the action of $S$ on $M$, the continuity of $f$,
and the continuity of $p$, there is an open neighborhood $W$ of $s$ in
the locally compact semigroup $S$ with compact closure $C=\overline{W}$
such that the second term $$p(f(x.t)-f(x.s))$$ is strictly less than
$\varepsilon/2$ for all $t\in W$ and $x\in K$.

Next, we note that the set $K.C$ is a compact subset of $M$, being the
image of $K\times C$ under the continuous action map $M\times S\to
M$. For all $g\in U_{p,K.C,\varepsilon/2}(f)$, we then have
$p(g(x.t)-f(x.t))<\varepsilon/2$ for all $t\in C$ and $x\in K$, by
definition. We conclude that
$p((t.g)(x)-\varphi(x))<\varepsilon/2+\varepsilon/2=\varepsilon$ for
all $t\in W$ and $g\in U_{p,K.C,\varepsilon/2}(f)$, which shows
\begin{equation*}
\Phi\big(W\times U_{p,K.C,\varepsilon/2}(f)\big)\subseteq
U_{p,K,\varepsilon}(\Phi(f,s)),
\end{equation*}
that is, continuity of $\Phi$.
\end{proof}

\begin{theorem}\label{thm:action}
Let $M$ be a complex Banach manifold, $V$ be a complete locally 
convex space, $G$ be a
\emph{finite dimensional} complex Lie group, and
$G \times M \to M$ a holomorphic right action. Then the 
action 
\begin{equation*}
G\times {\cal O}(M,V)\to {\cal O}(M,V),\quad (g.f)(x)=f(g^{-1}.x) 
\end{equation*}
is holomorphic with respect to the compact open topology on $\cO(M,V)$. 
\end{theorem}

\begin{proof} (\cite[Lemma~2.2.6]{Mue06}) 
This is a less general version of \cite[Theorem~III.14]{Ne01}, 
so the proof can be found there. However, the proof 
uses the defective Lemma~III.2(iii) from the same source, which we can
now replace with our Lemma~\ref{lem:action}. 
\end{proof}

\end{document}